\documentclass[reqno]{amsart}
\usepackage{amsfonts,amscd, amssymb,  epsf, epsfig, graphicx}

\def\noi{\noindent}

\newcommand{\bea}{\begin{eqnarray*}}
\newcommand{\eea}{\end{eqnarray*}}

\newtheorem{lem}{LEMMA}[section]
\newtheorem{theo}[lem]{THEOREM}
\newtheorem{claim}[lem]{Claim}
\newtheorem{coro}[lem]{COROLLARY}
\newtheorem{prop}[lem]{PROPOSITION}
\newtheorem{proposition}[lem]{PROPOSITION}
\newtheorem{definition}[lem]{DEFINITION}

\newtheorem{remark}[lem]{Remark} 
\newtheorem{question}[lem]{Question}
\newtheorem{example}[lem]{Example}
\newtheorem{conjecture}[lem]{Conjecture}
\newtheorem{problem}[lem]{Problem}
\def\CQFD{\hfill \vrule width 7pt height 7pt depth 1pt}

\newcommand{\ra}{\rightarrow}

\def\CC{{\rm\kern.24em\vrule width.02em height1.4ex depth-.05ex\kern-.26em C}}
\def\QQ{{\rm\kern.24em\vrule width.02em height1.4ex depth-.05ex\kern-.26em Q}}
\def\PP{{\rm\kern.24em\vrule width.02em height1.4ex depth-.05ex\kern-.26em P}}
\def\Rr{{\rm I\kern-.2em R}}
\def\ZZ{{\rm\kern.26em\vrule width.02em height0.5ex
depth0ex\kern.04em\vrule width.02em height1.47ex depth-1ex\kern-.34em Z}}
\def\BB{{\rm\kern.24em\vrule width.02em height1.4ex depth-.05ex\kern-.26em B}}
\def\RR{\hspace{.065in}\rm{\vrule width.02em height1.55ex
depth-.07ex\kern-.3165em R}}
\def\Ibb#1{{\rm I\kern-.23em#1}}

\def\noi{\noindent}

\def\noi{\noindent}

\newcommand {\cal} {\mathcal}
\newcommand {\be} {\begin{equation}}
\newcommand {\ee} {\end{equation}}

\newcommand {\beas} { \begin{eqnarray*}}
\newcommand {\eeas} {\end{eqnarray*}}

\begin{document}

\centerline{HARMONIC CURRENTS OF FINITE ENERGY AND LAMINATIONS}

\medskip

\centerline{by}

\medskip

\centerline{John Erik Forn\ae ss\footnote{The first author is supported by an NSF grant.
} and Nessim Sibony
}

\bigskip

%\vspace{.5cm}

%\centerline{         by }

%\vspace{.5cm}

%\centerline{ John Erik Forn\ae ss$^*$ and Nessim Sibony}

%\vspace{.5cm}

%\today

%\tableofcontents

%\vspace{.5cm}

Abstract: We introduce, on a complex K\"{a}hler manifold $(M,\omega),$
a notion of energy for harmonic currents of bidegree $(1,1).$
This allows us to define $\int T \wedge T \wedge \omega^{k-2},$
for positive harmonic currents. We then show that for a lamination
with singularities of a compact set in $\PP^2$ there is a unique
positive harmonic current which minimizes energy. If $X$ is a compact
laminated set in $\PP^2$ of class ${\mathcal C}^1$ it carries a
unique positive harmonic current 
$T$ of mass $1.$ The current $T$ can be obtained
by an Ahlfors type construction starting with a arbitrary leaf of $X.$

%\maketitle

%\authors{John Erik Forn\ae ss, Nessim Sibony}

%$^*$ The first author was supported by an NSF grant. 

\section{Introduction}

Let $X$ be a compact set in a complex manifold $M.$ If $X$ is laminated
by Riemann surfaces, a result due to L. Garnett implies the existence
of a positive current $T$ of bidimension $(1,1)$ which is harmonic
i.e. such that $i\partial \overline{\partial}T=0.$ Moreover in a flow box
$B,$ the current can be expressed as

$$
T=\int h_\alpha [V_\alpha]d\mu(\alpha),\;\;\;\;\;\;\;\;\;\;(1)
$$

\noindent the functions $h_\alpha$ are positive and harmonic on the 
local leaves $V_\alpha$, and $\mu$ is a measure on the transversal.
On the other hand there is a well known problem. Does there
exist a compact laminated set $X$ in $\PP^2$ which is not a compact
Riemann surface? See [CLS], [Gh]
and [Z] where the problem is dicussed. If such an $X$ does not exist
then the closure of any leaf $L$ of a holomorphic lamination in $\PP^2$
will contain a singularity. 

\bigskip

In this paper we study positive harmonic currents directed by a laminated
set with singularities. More precisely we consider compact sets $X$ laminated
by Riemann surfaces out of an exceptional set $E.$ We will assume
that $E$ is locally pluripolar and that $\overline{X \setminus E}=X.$
We will call such a set $(X,{\mathcal L},E)$ a laminated compact set
with singularities.
We consider on such sets harmonic currents $T$ of order $0$ and bidegree
$(1,1)$ in $\PP^2.$ The current $T$ can be written in the form

$$
T=c\omega+ \partial S + \overline{\partial S}
$$

\noindent where $\omega$ is the standard K\"{a}hler form on $\PP^2, c
\in \RR$, $S$ is a $(0,1)$ current. It turns out that $\overline{\partial}S$
depends only on $T$ and that for a positive closed current one can
define the energy $E(T)$ of $T$ by the following integral

$$
E(T)=\int \overline{\partial} S \wedge \partial \overline{S}
$$

\noindent and that $0 \leq E(T)<\infty.$ It is hence possible to introduce
a Hilbert space of classes of currents of finite energy. With this in hand
the integral

$$
Q(T)=\int T \wedge T
$$

\noindent makes sense for positive harmonic currents (not just the ones
associated to $(X,{\mathcal L},E)$) and has the
usual meaning when $T$ is smooth.
We then prove (Theorem 2.22):

\begin{theo}
Let $(X,{\mathcal L},E)$ be a laminated compact set with singularities.
There is a closed positive laminated current on $X$ or
there is a $\underline{unique}$ positive harmonic laminated current $T$
on $X$ minimizing energy.
\end{theo}

We then study the geometric intersection of laminated currents. We show,
see Proposition 2.16 and Theorem 5.2:

\begin{theo}
If $X$ is a ${\mathcal C}^1$ laminated compact set in $\PP^2$, then
$X$ carries a $\underline{unique}$ laminated positive harmonic current $T$. 
The class of the current
$T$ is extremal in the cone of positive harmonic currents. Moreover
$\int T \wedge T=0$.
\end{theo}

When $X$ is not ${\mathcal C}^1$ we have to assume "finite transverse
energy" to get the result. The main tool is that the quadratic form
$Q$ is negative definite on the hyperplane $\{T; \int T \wedge \omega=0\}.$

\medskip

The current can be obtained using a variant of Ahlfors' technique. Let
$\Phi:\Delta \rightarrow L$ be the universal covering map from the unit disc
to a leaf $L.$ Let $G_r(z)=\frac{1}{2\pi} \log^+ \frac{r}{|z|}.$
Define $T_r:= (\Phi)_*(G_r[\Delta]), r<1.$ If $A(r)$ is the mass of $T_r$
we get that 

$$
T=\lim 
_{r \rightarrow 1} \frac{T_r}{A(r)}.
$$

For that purpose we need to estimate the derivative of $\Phi$,
the estimates are valid for any laminated set.

The main question that is left open in our approach is to estimate

$$
c(\PP^2):= \inf \{ \int T \wedge T; \int T \wedge \omega=1,T \geq 0,
i\partial \overline{\partial}T=0\}.
$$
If $c(\PP^2)>0,$ then there is no ${\mathcal C}^1$ laminated set in $\PP^2.$

\section{Harmonic currents}

A subset $Y$ of a complex manifold $M$ is laminated by Riemann surfaces 
if it admits
an open covering $\{U_i\}$ and on each $U_i$ there is a homeomorphism
$\phi_i=(h_i,\lambda_i):U_i \rightarrow \Delta \times T_i$ where $\Delta$ is the
unit disc and $T_i$ is a topological space. The $\phi_i^{-1}$ are holomorphic
in $z.$ Moreover, 
$$\phi_{ij}(z,t)=\phi_j \circ \phi_i^{-1}(z,t)=(h_{ij}(z,t),\lambda_{ij}(t))$$
\noindent where the $h_{ij}(z,t)$ are holomorphic with respect to $z.$
When $T_i$ is in a Euclidean space and $\phi_{i}$ extend to ${\mathcal
C}^k$ diffeomorphisms, we say that the lamination is ${\mathcal C}^k.$
We call the $U_i$ flow boxes, $\{\lambda_i=t_0\}$ is a plaque.
A leaf is a minimal connected set such that if $L$ intersects a plaque,
then $L$ contains the plaque. We only consider the case when we have
a compact set $X$ contained in a K\"{a}hler manifold $M$ of dimension $k.$
We assume that $X$ contains a closed set $E$ such that $\overline{X \setminus E}=X$ and
$X \setminus E$ is a lamination ${\mathcal L}$ by Riemann surfaces.
We call such a triple $(X,{\mathcal L},E)$ a lamination with singularities,
$E$ is the singular set.
We will say that $(X,{\mathcal L},E)$
is oriented if there are continuous nonvanishing $(1,0)$ form
$\gamma_j,j=1,\dots,k-1$ defined on $X \setminus E$ such that $\gamma_j
 \wedge [\Delta_\alpha]=0$
for every plaque $\Delta_\alpha,$ $[\Delta_\alpha]$ denotes the current
of integration on the disc $\Delta_\alpha.$ 
We only consider sets $E$ which are pluripolar, more precisely,
for any $p\in E$, there are small balls $B(p,r)$ and $u,$ a plurisubharmonic
function on $B(p,r)$ such that $E \cap B(p,r)
\subset \{u=-\infty\}$ and $u$ is not identically $-\infty$
on $K \cap B(p,r).$ When $E \cap B(p,r)=\{u=-\infty\}$ we say
that $E$ is locally complete pluripolar. 
A positive  current $T$ of bidimension $(1,1)$ with support in $X$
is said to be directed by ${\mathcal L}$ if on any open
set $U$ where ${\mathcal L}$  is defined by nonvanishing
continuous $(1,0)$ forms $\gamma_j,j=1 \dots,k-1 $
i.e. $\gamma_j \wedge [\Delta_\alpha]=0,$
we have
$$ T \wedge \gamma_j=0.$$

To introduce the notion of a minimal set for $X$ we need the following 
assumption on the family of leaves $L.$ 

$$
 {\mbox{There is a neighborhood }}\; {\cal V}\; {\mbox{of E,
so that no leaf is contained}}\\
{\mbox{ in}}\; {\cal V}. \;\;\;\;\;\;\;\;(*)
$$

\begin{definition}
A minimal set for $(X,{\cal L}, E)$ is a compact subset $Y\subset X$
such that $Y$ is not contained in $E$,
moreover $Y\setminus E$ is a union of
leaves $L$ and for every leaf $L \subset Y\setminus E$, $\overline{L}=Y.$ 
\end{definition}

\begin{proposition}
Let $(X,{\cal L},E)$ satisfy assumption $(*).$ Then there are minimal sets
$Y \subset X.$ Two different minimal sets intersect only on $E.$ If M is a 
surface where the Levi problem is solvable and $E$ is locally contained
in a complex hypersurface, then any 
two minimal sets intersect.
\end{proposition}

{\bf Proof:}
Let ${\cal V}$ be an open neighborhood of $E$ such that no leaf is contained
in ${\cal V}.$ Let $(X_\alpha)$ be an ordered decreasing chain saturated for 
${\cal L}.$  Let $X'_\alpha:= X_\alpha \cap (X \setminus \cal V) .$
Then $X'_\alpha \neq \emptyset.$ Hence $\cap X_\alpha \neq \emptyset.$
So Zorn's Lemma applies. This shows that minimal sets exist.\\

It follows that each $M \setminus Y$
is locally pseudoconvex away from $E$, and since $E$ is locally contained
in a complex hypersurface,
$M\setminus Y$ is pseudoconvex [GR]. If the Levi problem is solvable
on $M,$ i.e. if pseudoconvex domains are Stein,
each $M\setminus Y$ is Stein. Since Stein manifolds cannot have two ends,
it follows that any two minimal sets must intersect. \\

\CQFD\\

\begin{remark}
Oka solved the Levi problem in $\PP^k.$ See [E] for a proof
of this and some generalizations. 
The condition on $E$ in the above Proposition, 
can be relaxed to assuming that $E$ is meager
in the sense of [GR].
\end{remark}

\begin{example}
Let ${\cal L}_\alpha$ be the foliation in $\PP^2$ defined by 
$wdz-\alpha zdw=0,$ with $\alpha$ irrational. Then $Y_c=
\overline{\{|z|=c|w|^\alpha\}}$ is minimal for every $c,$
the closure is in $\PP^2.$ The associated positive closed current
$T$ is defined by $\pi^* T=i\partial \overline{\partial} u, u(z,w,t)=
\log (\max \{|z||t|^{\alpha-1}, |w|^\alpha\})$ if $\alpha>1$ and is directed
by ${\mathcal L}_\alpha.$
\end{example}

\begin{definition}
Let $M$ be a complex manifold of dimension $k$. For $0\leq
p,q\leq k,$ let $T$ be a $(p,q)$ current of order $0.$
We say that $T$ is harmonic if $i \partial \overline{\partial} T=0.$
\end{definition}

\medskip

Observe that if $T$ is harmonic then $\overline{T}$, the
conjugate, is also harmonic. A current is real if $T=
\overline{T},$ in which case $p=q.$

\subsection{Decomposition of harmonic currents}
 We want to prove a representation theorem for real harmonic 
currents on compact K\"{a}hler manifolds.

\medskip

There is also a notion of $\Box-$harmonic forms in the $\overline{\partial}$
literature [FK]. These are smooth $(p,q)$ forms $\Omega$ for which
$\Box \Omega=(\overline{\partial} \overline{\partial}^*+
\overline{\partial}^*\overline{\partial}) \Omega=0.$
These forms consist of the common null space
of $\overline{\partial}$ and $\overline{\partial}^*$.
Note that when $p=q$, $\overline{\Box}=\Box$ so the conjugate of 
a $\Box-$harmonic
form is also $\Box-$harmonic. Since $\Box-$ harmonic forms are $\overline{\partial}$
closed they are also harmonic in the sense of currents as defined above.
Recall that on a compact K\"{a}hler manifold, for a closed current $u$ of
bidegree $(p,q),$ the following are equivalent:\\

\noindent i) $u$ is exact,\\
ii) $u$ is $\overline{\partial}$ exact,\\
iii) $u$ is $\partial \overline{\partial}$ exact.\\

See Demailly [De], p 41 for smooth forms. The proof is the same
for currents since cohomology groups for currents and smooth forms
are the same. This follows from the deRham Theorem and the fact that
for K\"{a}hler manifolds the Dolbeault cohomology groups
$\oplus_{p+q=k} H^{p,q}(X,\CC)$ and $H^k_{{\mbox{DR}}}$ are isomorphic
[De] p.42.

\begin{proposition}
Let $T$ be a harmonic $(p,q)$ current on a compact
K\"{a}hler manifold $M$ of dimension $k.$ Then
$$
T=\Omega+\partial S+\overline{\partial} {R}\;\;\;\;\;\;\;(2)
$$
\noindent where $\Omega$ is a unique 
closed smooth $\Box$-harmonic
form of bidegree $(p,q)$, and $S$ is a current
of bidegree 
$(p-1,q),$ $R$ is of bidegree $(p,q-1).$
When $dT$ is of order $0,$ we can choose $S,R$ depending linearly on $T$ and 
$$
L:T \rightarrow (\Omega,S,R)
$$
is continuous in the topology of currents. When $T$ is real we can choose
$R=\overline{S}.$
\end{proposition}

{\bf Proof:}
The current $\overline{\partial}T$ is $\partial-$
closed, hence $d-$ closed and is $\overline{\partial}$ exact. It follows
from Lemma 8.6 in [De] that $\overline{\partial}T$ is $\partial \overline{
\partial}$ exact so there is an $S_0$ of order $0$ and
of bidegree $(p-1,q)$
such that 
$$
\overline{\partial}T = \overline{\partial}\partial S_0.
$$

Let $\Omega$ be a smooth $\Box-$ harmonic representative
of the Dolbeault cohomology class of $T-\partial S_0$.
(See [De] and [V].)

Then $T-\partial S_0-\Omega $ is
$\overline{\partial}$ exact. Hence
$$
T=\Omega+\partial S_0+\overline{\partial}R.
$$

The choices of $S_0$ and $R$ can be made linearly since there is an
explicit inverse for $\partial \overline{\partial}$ and hence
$\overline{\partial}$, on currents cohomologous to zero. See
Dinh-Sibony [DS] where an explicit kernel is given.
If $T$ is real we obtain that $T$ can be expressed as claimed.

The current $T$ acts on $H^{n-p,n-q},$ the Dolbeault cohomology group,
because of the $\partial \overline{\partial}$ lemma. Hence $\Omega$
is uniquely determined by $T.$

\CQFD\\

\medskip

\begin{prop}
 If $T$ is as in $(2)$, then $S,R$ are not unique, but any other
$S',R'$ can be obtained as $S'=S+\omega'+\partial v+
\overline{\partial} u,$ similarly for $R'.$
\end{prop}

{\bf Proof:} The cohomology class $\Omega$ is defined
uniquely. If $S',R'$ is another solution 
we get $\partial(S-S')+\overline{\partial} (R-R')
=0$. 
Assume $\partial \sigma+\overline{\partial}\overline\sigma'=0.$
Then $\sigma$ is harmonic. Using the above construction for a harmonic
$(p-1,q)$ form, we get

$$
S-S'=\sigma=\omega'+\partial v+\overline{\partial} u.$$

Hence

$$
S=S'+(\omega'+\partial v+\overline{\partial} u ).$$

\CQFD\\

\begin{coro}
Let $T$ be a harmonic current of degree $(1,p)$ on $(M,\omega).$
Let $T=\Omega+\partial S+\overline{\partial} R$ be any decomposition
as in (1). Then $\overline{\partial}
S$ is uniquely determined by $T.$ If $p=1$ and $T$ is real, then
$T$ is closed if and only if $\overline{\partial} S=0.$
\end{coro}

{\bf Proof:}
If $S',R'$ also satisfy $T=\Omega+\partial S'+\overline{\partial} R'$
then $S'=S+\omega'+\overline{\partial} u$ for bidegree reasons. Consequently
$\overline{\partial}S'=\overline{\partial}S.$
Assume that $T$ is a real $(1,1)$ current. If $T$ is closed, then
$T=\Omega+i\partial \overline{\partial}u$, i.e. $\overline{\partial}S=0.$
Conversely, if $\overline{\partial}S=0$ and $T=\Omega+\partial S+
\overline{\partial S},$

$$
<T,\overline{\partial}\theta>=-<\overline{\partial}T,\theta>
=<-\overline{\partial}\partial S,\theta>=0.$$

Hence $T$ is closed.\\

\CQFD\\

\subsection{Energy of harmonic currents}

In this paragraph we introduce a notion of energy of harmonic currents of
bidegree $(1,1)$ on $\PP^k.$ On $\PP^k$ we consider the standard 
K\"{a}hler form $\omega$ normalized so that $\int \omega^k=1.$
Recall that $H^{p,q}=0$
except when $p=q$ in which case $H^{p,p}$ is generated by $\omega^p.$

\medskip

We showed above that if $T$ is a real harmonic $(1,1)$ current on $\PP^k,$ then
it can be represented as

$$
T=c\omega+\partial S+\overline{\partial}\overline{S}
$$

\noindent with $S$ of bidegree $(0,1)$. We have $c\in \RR$, 
$c=\int T \wedge \omega^{k-1}.$
We define the energy $E(T)=E(T,T)$ of $T$ as
$$
E(T,T)=\int \overline{\partial}S\wedge \partial \overline{S}\wedge \omega^{k-2}
$$
\noindent when $\overline{\partial}S \in L^2$.
 Observe that
$\overline{\partial}S$
is a $(0,2)$ form. Hence $0\leq E(T,T) < \infty.$ We have seen in
Corollary 2.8 that the energy depends on $T$ only. 

\medskip

We define ${\cal H}_e$ to be the space of real harmonic $(1,1)$ currents on
$\PP^k$ of 
finite energy. We consider on ${\cal H}_e$ the (real) inner product and 
semi norm

\bea
<T_1,T_2>_e & = & \left[ \int T_1\wedge \omega^{k-1}\right]
 \left[ \int T_2\wedge \omega^{k-1}\right]+\frac{1}{2}
\int \overline{\partial}S_1\wedge \partial \overline{S}_2 \wedge \omega^{k-2}\\
& + & \frac{1}{2}\int \overline{\partial}S_2\wedge \partial \overline{S}_1
 \wedge \omega^{k-2}\\
\|T\|^2_e & = & \left| \int T \wedge \omega^{k-1}\right|^2+
\int \overline{\partial} S \wedge \partial \overline{S}\wedge \omega^{k-2}\\
\eea

\begin{lem}
Let $T=c\omega+\partial S+\overline{\partial S},
\overline{\partial}S\in L^2,$ a $(1,1)$ real harmonic current
of order $0$ in $\PP^k.$ Then
$\|T\|_e=0$ if and only if  $T=i \partial \overline{\partial}u,$
for $u \in L^1,$ $u$ real.
\end{lem}

{\bf Proof:}
If $\|T\|_e=0$, then $\overline{\partial}S=0$ and $c=0.$ Then
$T$ is closed, 
hence $T$ is exact and therefore $T=i\partial \overline{\partial}u$
[De]. Regularity of the Laplace equation shows that $u\in L^1.$
Conversely, suppose that $T=i \partial 
\overline{\partial} u,$ $u\in L^1,$ $u$ real, so
we can set $S=
\frac{1}{2}i\overline{\partial} u.$ Then $\overline{\partial}S=0,$
hence $\int \overline{\partial}S \wedge
\partial \overline{S}=0.$ Clearly also, $\int T \wedge \omega=0$
so  $\|T\|_e=0.$\\

\CQFD\\

\medskip

\begin{lem}
There is a constant $C=C_k$ so that 
if $T$ is a real harmonic $(1,1)$ current of order $0$ on $\PP^k, k \geq 2,$
with finite energy,
then there is an element $\tilde{T}$ in the equivalence class of 
$T$, i.e. $\|T-\tilde{T}\|_e=0,$ which can be written
as $\tilde{T}=c\omega+\partial S+ \overline{\partial S}$
with $S,\partial S, \overline{\partial}S \in L^2,$
$\|S\|_{L^2}, \|\partial S\|_{L^2}, \|\overline{\partial} S\|_{L^2}
\leq C \|T\|_e.$ Hence $T=\tilde{T}+i\partial \overline{\partial}u$
and $i\partial \overline{\partial }u$ is of order $0.$
\end{lem}

{\bf Proof:} We can write $T=c\omega+ \partial S_1+ \overline{ \partial S}_1$
with $\overline{\partial}S_1\in L^2.$ 
Since $ \overline{\partial} S_1 $ is in $ L^2$, we can find an
$S\in L^2_{01}$ for which 
$\overline{\partial} S_1=\overline{\partial} S.$
Moreover, $\partial S\in L^2$ as well since, by Hodge theory
we gain one derivative
by solving the $\overline{\partial}$ problem 
$\overline{\partial} S=\overline{\partial} S_1.$

Since $H^{(0,1)}(\PP^k)=0,$ there is a distribution $v$ for which 
$\overline{\partial}v= S_1-S.$ 
Therefore we have the decomposition

\bea
T & = &  c\omega+[\partial(S+\overline{\partial}v)]+\overline{
[\partial(S+\overline{\partial}v)]}\\
& = & c\omega+\partial S+\overline{\partial S}+\partial
\overline{\partial}v+\overline{\partial}\partial \overline{v}\\
& = & \tilde{T}+i \partial \overline{\partial} \left( \frac{
v-\overline{v}}{i} \right)\\
\eea

The distribution  $u:= \frac{v-\overline{v}}{i}$ is real.
Since $T,\tilde{T}$ have order $0$, $i \partial \overline{\partial}u$
is also of order $0.$

Since $\|\partial{\overline{\partial}}u\|_e=0,$ it follows that
$\tilde{T}:=  c\omega+\partial S+\overline{\partial S}$ is in the equivalence
class of $T$ as desired.

The $L^2$ estimates are classical [FK].\\

\CQFD\\

Let $H_e$ denote the quotient space of equivalence classes $[T]$
in ${\mathcal H}_e.$

\begin{proposition} The space $H_e$ is a 
real Hilbert space. Every element $[T]$ in $H_e$ can be represented
as
$$
T=c\omega+\partial S+\overline{\partial} \overline{S}
$$

\noindent where $S$ is a $(0,1)$ form in $L^2,$ with $\partial S$ and
$\overline{\partial}S$
in $L^2.$ Convergence in $H_e$ implies weak convergence of currents:
If $[T_n] \rightarrow 0$ in $H_e$ then there are representatives $\tilde{T}_n
\in [T_n]$ such that $\tilde{T}_n \rightarrow 0$ 
in the weak topology of currents. In fact the mass norms $\|\tilde{T}_n\|
\rightarrow 0.$
\end{proposition}

{\bf Proof:}
We show first that $H_e$ is complete. Let $\{[T_n]\}$ be a Cauchy sequence
of equivalence classes, $\lim_{n,m \ra \infty} \|T_n-T_m\|_e =0.$
We can suppose $\|T_{n+1}-T_n\|_e<\frac{1}{2^n}.$ Inductively, we can 
(Lemma 2.10) choose 
representatives $\tilde{T}_n$ so that 

$$
\tilde{T}_{n+1}=\tilde{T}_n+c_n \omega+\partial S_n+ \overline{\partial S_n},
|c_n|, \|\partial S_n\|_{L^2}, \|S_n\|_{L^2}, 
\|\overline{\partial}S_n\|_{L^2}
\leq C \frac{1}{2^n}.
$$

Hence $\{[T_n]\}$ converges in $H_e.$ This shows that $H_e$ is complete.
The last statement is similar.\\

\CQFD\\

We will next introduce a notion of wedge product of real
harmonic currents. Let $T,T'$ be representatives of equivalence classes
as above,

$$
T=c\omega+\partial S+ \overline{\partial S},
T'=c'\omega+\partial S'+ \overline{\partial S'}.
$$

Then a formal calculation gives:

\bea
\int T \wedge T'\wedge \omega^{k-2} & = &
cc'\int \omega^k+
\int \partial S \wedge \overline{\partial S'}\wedge \omega^{k-2}+
\int \overline{\partial S} \wedge {\partial S'}\wedge \omega^{k-2}\\
& = & <T,\omega^{k-1}><T',\omega^{k-1}>-
\int \overline{\partial} S \wedge \partial \overline{S'}\wedge \omega^{k-2}\\
& - & 
\int \partial \overline{S} \wedge \overline{\partial} S'\wedge \omega^{k-2}\\
\eea

Notice that if $T,T'$ have finite energy, the last expression
is well-defined. We define in this case the quadratic form $Q(T,T')$
for currents $T,T'$ of finite energy:

$$
Q(T,T')=
<T,\omega^{k-1}><T',\omega^{k-1}>-
\int \overline{\partial} S \wedge \partial \overline{S'}\wedge \omega^{k-2}-
\int \partial \overline{S} \wedge \overline{\partial} S'\wedge \omega^{k-2}
$$

\noindent and motivated by the formal calculation, we define

$$
\int T \wedge T'\wedge \omega^{k-2}:=Q(T,T')
$$

\noindent when $T,T'$ are harmonic $(1,1)$ current on $\PP^k$
with finite energy.
Recalling the definition of energy,
we get

$$
\int T \wedge T\wedge \omega^{k-2}=Q(T,T)= <T,\omega^{k-1}>^2
-2E(T,T).
$$

Note that $Q(T,T')$ is well defined on equivalence classes in $H_e.$

\begin{theo}
Any positive harmonic current $T$ of bidegree $(1,1)$ on $\PP^k$ is
of finite energy.  Also $\int T \wedge T \wedge \omega^{k-2}\geq 0$
for these currents.
The quadratic form $Q(T,T')$
is continuous on $H_e.$
\end{theo}

{\bf Proof:} Assume first that $T$ is smooth, 
$T=c\omega+\partial S
+\overline{\partial}\overline{S}$ with $c \geq 0$
and $S$ smooth. We get after integration by parts:

\bea
\int T \wedge T \wedge \omega^{k-2} & = & 
c^2\int \omega^k+2\int \partial S \wedge \overline{\partial S}\wedge 
\omega^{k-2}\\
& = & c^2 \int \omega^k-2\int \overline{\partial} S \wedge \partial 
\overline{S} \wedge \omega^{k-2}\\
& \geq & 0\\
\eea

So
$$
2\int \overline{\partial} S \wedge \partial \overline{S}\wedge \omega^{k-2}
\leq |<T,\omega^{k-1}>|^2.
$$

In general, we can still write $T=c\omega+\partial S
+ \overline{\partial} \overline{S}$ by Proposition 2.6.
Let $S_\epsilon$ be a regularization of $S$ and define
$T_\epsilon=c\omega+\partial S_\epsilon+\overline{\partial}
\overline{S}_\epsilon.$ Here we use the classical regularization
for a current $S$, $$S_\epsilon= \int_{U(k+1)}\rho_\epsilon(g)g_*Sd\nu(g),$$
\noindent where $\nu$ is the Haar measure on $U(k+1)$ acting as automorphisms
of $\PP^k$ and $
\rho_\epsilon$ is an approximation of unity in $U(k+1).$

Then $T_\epsilon$ is still positive and harmonic. We get that
$$
2\int \overline{\partial} S_\epsilon \wedge \partial \overline{S}_\epsilon
\wedge \omega^{k-2} \leq |<T,\omega^{k-1}>|^2\;\;\;\;\;\; (3)
$$

Since $S_\epsilon \rightarrow S$ weakly,
 $\overline{\partial}S_\epsilon
\rightarrow \overline{\partial}S$ converges weakly in $L^2$, because
$\{{\overline{\partial}S}_\epsilon\}_\epsilon$ is bounded in $L^2.$
Hence $\overline{\partial}S\in L^2$ and
$$
\int \overline{\partial}S \wedge \partial \overline{S}\wedge \omega^{k-2}
\leq \underline{\lim}\int \overline{\partial}S_\epsilon
\wedge \partial \overline{S}_\epsilon\wedge \omega^{k-2}
\leq 1/2 |<T,\omega^{k-1}>|^2.$$

\noindent Hence $T$ has finite energy and $\int T \wedge T \wedge \omega^{k-2}\ge 0.$

If $T=c\omega+\partial S+\overline{\partial} \overline{S}$ and
$T'=c\omega+\partial S'+\overline{\partial} \overline{S'}$ then
$Q(T,T')=cc'\int \omega^k-
\int \overline{\partial}S \wedge \partial \overline{S'}\wedge \omega^{k-2}-
\int \overline{\partial}S' \wedge \partial \overline{S}\wedge \omega^{k-2}.$
It is clear that $|Q(T,T')| \leq 2 \|T\|_e\|T'\|_e.$ Hence $Q$
is continuous on $H_e.$\\

\CQFD\\

\begin{remark}
Since our regularization is a convolution, we also have $\overline{\partial}
S_\epsilon \rightarrow \overline{\partial}S$ in $L^2$ and hence
$E(T_\epsilon,T_\epsilon)\rightarrow E(T,T).$
\end{remark}

\begin{coro}
On $H_e$, the quadratic form
$$
Q(T_1,T_2)=\int T_1 \wedge T_2 \wedge \omega^{k-2} 
$$
is strictly negative definite on the hyperplane
${\cal H}=\{T; \int T \wedge \omega^{k-1}=0\}.$ Consequently, if $T,T'$
are positive harmonic currents, non proportional, then
$\int T \wedge T'\wedge \omega^{k-2}>0.$ If
$\int T \wedge \omega^{k-1}=1$ and $T \geq 0$ then $T$ is 
non closed if and only
if $0 \leq \int T \wedge T \wedge \omega^{k-2}<1.$  
\end{coro}

{\bf Proof:}
If $\int T \wedge \omega^{k-1}=0,$ then
$$
Q(T,T):=
-2\int \overline{\partial}S \wedge \partial \overline{S} \wedge \omega^{k-2}
\leq 0,$$ it is zero only if $T=0$ in $H_e.$ Hence $Q$ is strictly
negative definite on ${\mathcal H}.$ 
Suppose the space generated by $T',T$ is of dimension $2.$ There is an $a>0$,
$\int (T'-aT)\wedge \omega^{k-1}=0. $
Hence 
\bea
0 & > & Q(T'-aT,T'-aT)\\
& = & Q(T',T')+a^2 Q(T,T)-2aQ(T',T)\\
\eea

\noindent Since $Q(T',T'), Q(T,T) \geq 0$ by Theorem 2.12, it
follows that $Q(T',T)>0.$ 

The last part is an immediate consequence of Corollary 2.8.\\

\CQFD\\

\begin{prop}
The function $T \rightarrow Q(T,T)$ is upper semi continuous
in the weak topology on positive harmonic currents and is strictly
concave on $\{T \wedge \omega=1\}.$
\end{prop}

{\bf Proof:}
If $T_n \rightarrow T$ weakly, we
have seen as above in the proof of Theorem 2.12 
that $\int \overline{\partial}S \wedge \partial \overline{
S} \leq \underline{\lim}\int \overline{\partial}S_n
 \wedge \partial \overline{S}_n$, so $Q$ is upper semi continuous.
Concavity is clear because if $\int (T-T')\wedge \omega=0$,
$Q(T-T',T-T')<0,$ so $2Q(T,T')>Q(T,T)+Q(T',T').$ Hence
$Q(\frac{T+T'}{2}, \frac{T+T'}{2})> \frac{1}{2} Q(T,T)+
\frac{1}{2} Q(T',T').$ \\

\CQFD\\
 
\begin{prop}
In $\PP^k, k \geq 2,$
positive harmonic currents $T$ of bidegree $(1,1)$ satisfying
$Q(T,T)=0$ are extremal on $H_e$ in the cone of positive harmonic currents.
\end{prop}

{\bf Proof:}
Assume that $0\leq T'\leq cT$ for some $c>0.$

\bea
\int T' \wedge T \wedge \omega^{k-2} & = &
\lim \int T_\epsilon' \wedge T \wedge \omega^{k-2}\\
& \leq & \lim \int cT_\epsilon \wedge T \wedge \omega^{k-2}\\
& = & c \int T \wedge T \wedge \omega^{k-2}\\
& = & 0\\
\eea

Hence $[T']$ is proportional to $[T]$ by Corollary 2.14.\\

\CQFD\\

\begin{prop}
The map $T \rightarrow \overline{\partial}S$ is not continuous for the
weak topology on positive harmonic currents $T$ of bidegree $(1,1)$
and $L^2$ topology on $\overline{\partial}S.$ 
\end{prop}

{\bf Proof:}
Let $T=\omega+\epsilon(\partial S+\overline{\partial S})$ for a smooth
$(0,1)$ form $S$ supported in the unit bidisc. For $\epsilon>0$ small
enough, $T$ is positive and
$\int T \wedge T= \int \omega \wedge \omega -2\epsilon^2 
\int \overline{\partial}S
\wedge \partial {\overline{S}}<1$.
If the map $T\rightarrow \overline{\partial}S$
with weak topology on $T$ and $L^2$ topology on $\overline{\partial}S$
were continuous then if $T_n \rightarrow T_0$,
$\int T_n \wedge T_n \rightarrow \int T_0 \wedge T_0.$
Let $f$ be an endomorphism of $\PP^2$ of algebraic degree $d.$
Then
$f^*:  {{\cal H}_e \rightarrow \cal H}_e$ is a linear map of norm $d$
if the algebraic degree of $f$ is $d.$ Indeed
\bea
|\int f^*T \wedge \omega|^2 & = & |\int T\wedge f_*\omega|^2\\
& = & d^2|\int T \wedge \omega|^2\\
\eea
because $f_*\omega \sim d\omega.$ We also have

$$
\left|\int f^*(\overline{\partial}S \wedge \partial \overline{S})\right|
=d^2 \left| \int \overline{\partial}S \wedge \partial 
\overline{S}\right|
$$

\noindent This can be obtained by smoothing.

Therefore
$E(f^*T/d)=E(T)$ so
$\int f^*T/d \wedge f^*T/d=\int T \wedge T.$

Let $f[z:w:t]=[z^2:w^2:t^2]$ and $T=\omega+\epsilon(\partial S+
\overline{\partial S})$ as above. Then

\bea
1 & > & \int T \wedge T \\
& = & \int (f^n)^*T/2^n \wedge (f^n)^*T/2^n \\
& = & \int (f^n)^*\omega/2^n \wedge
(f^n)^*\omega/2^n-2\int \overline{\partial} (f^n)^*(S)/2^n
\wedge {\partial} (f^n)^* (\overline{S})/2^n.\\
\eea

\noindent If we choose $S=a(|z|,|w|)d\overline{z}$
it is easy to check that
$(f^n)^*S/2^n
\rightarrow 0$ weakly. Hence $(f^n)^*T/2^n$ converges
weakly to a closed current $A=\lim \frac{(f^n)^* \omega}{2^n}$
whose class in $H^{(1,1)}$ is $\omega.$
If $T \rightarrow \overline{\partial} S$ were continuous
the second integral would converge to zero.
Since $\int \omega \wedge \omega =1,$ we get 
that
the map $T \rightarrow \overline{\partial} S$ is not
continuous. \\

\CQFD\\

This is the justification for introducing the norm
on finite energy currents, which gives a different topology than
weak topology.

\medskip

L. Garnett has shown in [G] the existence of positive currents
$T$, satisfying $i \partial \overline{\partial} T=0$ and directed
by foliations. In [BS] a version of this result is given allowing
leaves to intersect. Here we are interested in constructing laminar
currents for a foliation with singularities that are only holomorphic
motions in flow boxes, the holomorphic case in treated in
[BS].

\begin{theo}
Let $(X,{\mathcal L},E)$ be a directed set with singularities
in $\PP^2.$ Then there is a laminated harmonic positive current
$T$ of the form $T=\int_\alpha h_\alpha [\Delta_\alpha]d\mu(\alpha)$
in flow boxes. Here $\mu(\alpha)$ is a measure on transversals,
$h_\alpha$ are strictly positive harmonic functions,
uniformly bounded above and below by strictly positive constants,
$h_\alpha$ are Borel measurable with respect to $\alpha.$
\end{theo}

{\bf Proof:} Let $\gamma$ be a continuous $(0,1)$ form such that
$\gamma \wedge [\Delta]=0$ for every plaque $\Delta.$ In [BS]
Theorem 1.4 a current $T \geq 0$ supported on $X$
satisfying $i \partial \overline{\partial} T=0$ and
$T \wedge \gamma =0$ is constructed. It is shown that the current
is laminar in flow boxes when the foliation is holomorphic.
We consider here the general case.

The Theorem follows from the next Lemma.

\begin{lem}
A positive harmonic current, directed by a holomorphic motion in a 
polydisc is a laminar current.
\end{lem}

{\bf Proof of the Lemma:} Let $B$ be a flow box.
In $B$, $T=\|T\|i\gamma \wedge \overline{\gamma}$
where $\|T\|$ is a positive measure. Assume $z=0$ is a transversal 
and $\pi$ is the projection along leaves in $B,$ on that transversal.
Let $(\nu_w)$ be the desintegration of $\|T\|$ along $\pi.$ 
If $\phi$ is a test form of bidegree $(1,1)$ supported on $B$,

$$
<T,\phi>=\int <\nu_w i\gamma \wedge \overline{\gamma},\phi> d\mu(w)
$$

\noindent with $\mu=\pi_*(\|T\|).$
This applies when $\phi=i \partial \overline{\partial} f$ and when
$f$ has support on an arbitrary open set of plaques, we have:
$<T,i \partial \overline{\partial}f>=0.$ Hence

\bea
0 & = & \int <\nu_wi\gamma \wedge \overline{\gamma},i\partial
\overline{\partial} f>d\mu(w)\\
& = & \int <\nu'(w)[D_w],\Delta_w f>
d\mu(w).\\
\eea

\noindent where $\Delta_w f$ denotes the Laplacian along the leaf through
$w.$ This extends to test functions which are continuous and
${\mathcal C}^2$ along leaves.
Hence $<\nu'(w)[D_w],\Delta_w f>=0$
$\mu$ a.e., so $\nu'(w)$ is a harmonic function on the leaf
$[D_w].$  \\

\CQFD\\

\begin{remark}
It follows from a Theorem by Skoda ${\mbox{[Sk]}}$ that no positive
harmonic $(1,1)$ current can have mass on a set of $2-$
dimensional Hausdorff measure $\Lambda_2=0.$
\end{remark}

\begin{theo} Let $(X,{\cal L},E)$ be a laminated set with singularities
in $\PP^2.$ There is a $\underline{\mbox{unique}}$ equivalence class
of harmonic currents 
directed by ${\cal L}$ of mass one and minimal energy.
\end{theo}

{\bf Proof:} Let $C_1=\{T; T \geq 0, \int T \wedge \omega=1,
T\; {\cal L}-\mbox{directed, harmonic}\}.$ Then $C_1$ is compact
in the weak topology of currents. From Theorem 2.18 we know that $C_1$
is nonempty.
 
The energy is a lower semi continuous function on $C_1$ by Proposition 2.15.
Since $C_1$ is compact  it now follows that $E(T,T)$ 
takes on a minimum value $c$ on $C_1$

If $E(T,T)=c$ and $\int T \wedge \omega=1$, then
$Q(T,T)=1-2c.$
If $T,T'$ are two elements of non colinear equivalence classes
of currents where the minimum $c$
is reached, then $Q(\frac{T+T'}{2}, \frac{T+T'}{2})>1-2c$
by strict concavity of $Q$ (Proposition 2.15),
a contradiction. So the minimum is unique.\\

\CQFD\\

We next show that under mild extra hypotheses, minimal equivalence classes
contain only one current.
Recall that a current on a laminated compact $X$ 
which in local flow boxes has the form
$\int h_\alpha [V_\alpha]d\mu(\alpha)$ for a positive measure
$\mu(\alpha)$ on the space of plaques, and $h_\alpha>0,$ harmonic
functions on plaques $V_\alpha$ is said to be a laminated positive harmonic
current. The current is closed and laminated if the $h_\alpha$
are constant.

\begin{theo}
Let $(X,{\mathcal L},E)$ be a directed set with singularities.
Suppose $E$ is locally complete pluripolar with $\Lambda_2(E)=0$
in $\PP^2$. Assume there is no nonzero
positive closed laminated current on $X.$ Consider the convex compact
set $C$ of laminated positive harmonic $(1,1)$ currents of mass $1$. Then 
there is a unique element $T$ in $C$ minimizing energy. The current $T$
is extremal in $C.$
\end{theo}

{\bf Proof:}
We know that $C$ is nonempty. Let $T_1,T_2 \in C$ be two minimizing currents.
We know by remark 2.20 that $T_j$ has no mass on $E.$
By Theorem 2.21, $[T_1]=[T_2]$ so 
$T_1-T_2=i\partial\overline{\partial}u$, hence $T_1-T_2$ is closed.

In a flow box we have $T_j=\int h_j^\alpha[V_\alpha]d\mu_j(\alpha),
j=1,2.$
Let $\nu(\alpha)=\mu_1+\mu_2,$ so $\mu_j=r_j(\alpha)\nu$. Then

$$
T_1-T_2= \int \left(h_1^\alpha r_1(\alpha)-h_2^\alpha r_2(\alpha)\right)
[V_\alpha]d\nu(\alpha).$$

Since $d(T_1-T_2)=0$ it follows that $$h_1^\alpha r_1(\alpha)-
h_2^\alpha r_2(\alpha) \equiv c(\alpha)$$ 

We decompose the measure $c(\alpha)\nu(\alpha)$ 
on the space of plaques, $c(\alpha)
\nu(\alpha)=\lambda_1-\lambda_2$ for positive mutually singular
measures $\lambda_j.$ Then

$$T_1-T_2= \int [V_\alpha]\lambda_1(\alpha)-\int [V_\alpha]\lambda_2(\alpha)
=T^+-T^-$$ for positive closed currents $T^\pm.$ These locally defined
currents fit together to global positive closed currents on $K\setminus E.$
Observe that the mass of $T^\pm$ is bounded by the mass of $T_1+T_2$.

Since $E$ is
locally complete pluripolar the trivial extensions of $T^\pm$ are also closed. 

 Consequently $T^\pm \equiv 0$ and $T_1=T_2.$
The fact that $T$ is extremal follows from the strict concavity of $Q$.\\

\CQFD\\

\begin{coro}
Let $(X,{\cal L},E)$ be a laminated set with singularities in $\PP^2$ with
$\Lambda_2(E)=0.$ One of the
following statements holds:\\
i) There is a closed positive laminated 
current of mass $1$ on $X.$\\
ii) There is a unique positive laminated harmonic current of mass one
on $X.$\\
iii) There is a positive laminated harmonic current $T$ of mass one on $X$
such that $\int T \wedge T >0.$ In particular the current $T_0$ 
minimizing energy satisfies $\int T_0 \wedge T_0>0.$
\end{coro}

{\bf Proof:} Assume i) and ii) fail. Then there are two
positive harmonic non closed currents $T_1,T_2$ of mass one which
are not colinear. We can assume that $\int (T_1-T_2)\wedge \omega =0.$ 
By Corollary 2.14,
$Q(T_1-T_2,T_1-T_2)<0.$ If $Q(T_1,T_1)$ or $Q(T_2,T_2)$ is strictly
positive, we are done since $\int T \wedge T= Q(T,T).$ If not, then
$Q(T_1,T_2)=\int T_1 \wedge T_2 >0.$ But then if $T:=\frac{T_1+T_2}{2},$
$Q(T,T)>0.$\\

\CQFD\\

\section{Intersections of laminar currents}

\subsection{${\cal C}^1$ laminations}

Here we assume that $X\subset \mathbb P^2$ is a laminated compact
covered by finitely many flow boxes $B_i.$
We suppose that $X$ locally extends to a ${\cal C}^1$ lamination of
an open neighborhood.
We can assume that $p:=[0:0:1]\in B_1 \subset X$ and that the leaf $L$ through
$p$ has the form $w={\cal O}(z^2)$. So we will assume that the lamination of $X$ is of the
form $w=w_0+f_{w_0}(z), f_{w_0}(0)=0,$ where the map $\Psi(z,w_0)=(z,
w_0+f_{w_0}(z))$ is a ${\cal C}^1$ diffeomorphism in a neighborhood
of $p.$

\medskip

Let $\Phi_\epsilon([z:w:t])=[z,w+\epsilon z:t]$ denote
a family of automorphisms of $\mathbb P^2.$ Notice that each
of these automorphisms fixes the $w$ axis.

For two graphs $L_1,L_2$ given by $w=g_1(z),w=g_2(z), z\in S
,$ we define the vertical distances over $S$ between the two
as $d_S^{max}(L_1,L_2)=\sup_{z\in S}|f_1(z)-f_2(z)|$
and $d_S^{min}(L_1,L_2)=\inf_{z\in S}|f_1(z)-f_2(z)|.$

\begin{theo}
There exists an integer $N$ independent of $\epsilon $
so that in any of the flow boxes $B_i$
local leaves $L_1$ and $\Phi_\epsilon(L_2)$ can at most
intersect in $N$ points, counted with multiplicity. Moreover there exist
neighborhoods $U_\epsilon$ of ${\mbox{Id}}$ in $U(3)$ so that the
same conclusion holds for $\Psi_1(L_1)$ and $\Psi_2(\Phi_\epsilon(L_2))$,
$\Psi_1, \Psi_2\in U_\epsilon.$  
\end{theo}

{\bf Proof:} Fix a $\delta>0.$ Let $L_{w}$ denote the leaf
through $[0:w:1]$ and let $L^\epsilon_w$ denote
its image under $\Phi_\epsilon.$ 
Say $L_{w_0}$ is given by $w=w_0+f_{w_0}(z)$ 
and $L^\epsilon_{w_0}$ is given by
$w=w_0+f_{w_0}(z)+\epsilon z.$ Note that the vertical
distance $d_S^{\mbox{max}}$ between  $L_{w_0}$ and $L^\epsilon_{w_0}$
is $|\epsilon| \delta$ at the boundary $S$ of the disc $|z|\leq 
\delta.$
Because the lamination is of class ${\mathcal C}^1,$ there
exists a constant $C>1$ so that if $(0,w_0),(0,w_1)\in B_1,$
then 

$$
\sup_{|z|\leq \delta}
 |w_0+f_{w_0}(z)-w_1-f_{w_1}(z)| \leq C
\inf_{|z|\leq \delta}
 |w_0+f_{w_0}(z)-w_1-f_{w_1}(z)|$$

If $L_1$ and $\Phi_\epsilon(L_2)$ intersect in a flow box $B_i$,
then $L_1$ and
$L_2$  must be at most $a|\epsilon|$ apart in $B_i$ for some
fixed $a.$ Let $c>0$ be any small constant. If
$L_1$,  $\Phi_\epsilon(L_2)$ intersect in $N$ points
in $B_i$ and $N$ is sufficiently large then
$L_1$ and  $\Phi_\epsilon(L_2)$ can be at most a distance
$c |\epsilon|$ apart from each other in $B_i.$ Note that the
same conclusion holds if the number of intersections is counted for
$\Psi_1(L_1)$ and $\Psi_2(\Phi_\epsilon(L_2))$ for a small enough
$U_\epsilon.$  
However, there is a
path of at most a fixed length along these leaves ending
in the flow box containing $[0:0:1]$.
It follows that continuing these leaves to this
neighborhood, they will have to stay at most $bc |\epsilon|$ apart,
for a fixed constant $b.$ Choosing $c$ small enough, we get
$bc<\frac{\delta}{2C}$. Let $L_1=\{w=w_1+f_{w_1}(z)\},$ 
$L_2=\{w=w_2+f_{w_2}(z)\}.$
Then $|w_2+f_{w_2}(z)+\epsilon z-w_1-f_{w_1}(z) |
<\frac{|\epsilon|\delta}{2C}$ when $|z| \leq \delta.$ Hence

\bea
d_{|z| \leq \delta}^{\mbox{max}}(L_1,L_2) & \leq &
C d_{|z| \leq \delta}^{\mbox{min}}(L_1,L_2)\\
 & \leq & C d_{z=0}^{\mbox{min}}(L_1,L_2)\\
& = & C d_{z=0}^{\mbox{max}}(L_1,\Phi_\epsilon (L_2))\\
& \leq & C d_{|z| \leq \delta}^{\mbox{max}}(L_1,\Phi_\epsilon (L_2))\\
& \leq & \frac{|\epsilon|\delta}{2}.\\
\eea

Applying this estimate when $|z|=\delta,$ we get 

\bea
\frac{|\epsilon|\delta}{2} & > &  \frac{|\epsilon|\delta}{2C}\\
& > & |w_2+f_{w_2}(z) +\epsilon z-w_1-f_{w_1}(z)|\\
& \geq & |\epsilon| \delta -|w_2+f_{w_2}(z)-w_1-f_{w_1}(z)|\\
& \geq & |\epsilon|\delta -\frac{|\epsilon|\delta}{2}\\
& = & \frac{|\epsilon|\delta}{2},\\
\eea

\noindent a contradiction.\\

\CQFD\\

\subsection{Laminations by holomorphic motions}

Now we consider the case of laminations which are not ${\mathcal C}^1.$
We recall the following result
by Bers-Royden [BR].

\begin{proposition}
We are given a lamination of a neighborhood of the unit polydisc
in $\CC^2.$ Assume that the leaves are of the following
form:

$$
L_t, t\in \CC, |t|<C, w=F_t(z), F_t(0)=t,F_0(z)\equiv 0.
$$

The map $\Phi(z)(t)=F_t(z)$ is a holomorphic motion and we have
the estimate:

$$
\frac{1}{K}|t-s|^{\frac{1+|z|}{1-|z|}} \leq |F_t(z)-F_s(z)|
\leq {K}|t-s|^{\frac{1-|z|}{1+|z|}}.
$$
\end{proposition}

\begin{theo}
Let $X\subset \mathbb P^2$ be a compact subset laminated
by Riemann surfaces. Then there exists a holomorphic family
$\Phi_\epsilon:\mathbb P^2 \rightarrow \mathbb P^2$ for
$\epsilon\in \mathbb C, \Phi_0\equiv \; {\mbox{Id}}$ with the
following properties. There are finitely many flow boxes
$\{B_i\}_{i=1,\dots,\ell}$ covering $X$ and an $\epsilon_0>0$ and a constant
$A$ such that if $L,L'$ are any local leaves in any flow box
$B_i$ then:\\
If $0<|\epsilon|<\epsilon_0$, the number of intersection points counted
with multiplicity of $L,\Phi_\epsilon(L')$ is at 
most $A \log \frac{1}{|\epsilon|}.$ Moreover there exist
neighborhoods $U_\epsilon$ of ${\mbox{Id}}$ in $U(3)$ so that the
same conclusion holds for $\Psi_1(L_1)$ and $\Psi_2(\Phi_\epsilon(L_2))$,
with $\Psi_1,\Psi_2 \in U_\epsilon.$
\end{theo}

{\bf Proof:}
We first choose a finite cover by flow boxes, $B_i$.
We can do this so that for each flow box there is a
linear change of coordinates in $\mathbb P^2$ so
that $[z:w:t]=[0:0:1]\in B_i \cap K$. Moreover, we can arrange
that if $L$ is any local leaf intersecting $\Delta(0,2)$
then $L\cap  \Delta(0,2)$ is contained in a local leaf
$\tilde{L}$ of the form $\{w=f_{\alpha}(z), |z|<3\},
(0,\alpha)\in \tilde{L},$
and $\|f_\alpha\|_\infty<3.$ Moreover we can assume that
each $\|f'_\alpha\|<.1$ and that $f'_0(0)=0.$ Redefining the flow boxes, we
can let
$B_i$ denote the union of those graphs over $|z|<3$
intersecting $\Delta(0,2).$ We can assume that
the smaller flow boxes $B_i'$ consisting of those
graphs over $|z|<1$ for which the graph is in $\Delta^2(0,1)$
already cover $K.$\\

Next we fix the coordinates $z,w,t$ on $\mathbb P^2$ used for
the first flow box $B_1$. Define the family $\Phi_\epsilon$ by
$$
\Phi_\epsilon [z:w:t]= [z:w+\epsilon z:t].
$$

\begin{lem}
There exists a $\; \delta>0$  and $C>0$ so that if
$w=f_\alpha (z), w=f_\beta (z)$ are two local leaves in $B_1$, then
$$
\frac{|\alpha-\beta|^2}{C}
 \leq |f_\alpha (z)-f_\beta (z)| \leq C|\alpha-\beta|^{1/2}, 
\; \forall \; z, |z| \leq \delta.
$$
\end{lem}

{\bf Proof:} This is a special case of the Bers-Royden result.\\

\CQFD\\

\begin{lem}
Let  $\epsilon_0>0$ be small enough. Then if
$L_1,L_2$ are leaves in the first flow box then
$d^{max}_{\{|z|\leq \delta\}}(\Phi_\epsilon (L_1),L_2)\geq |\epsilon|^3$
for all $|\epsilon|\leq \epsilon_0.$
\end{lem}

{\bf Proof:}
Let $L_i$ be given by $w=f_i(z), f_i(0)=w_i.$ Then
 $\Phi_\epsilon(L_1)$ is the graph $w=f_1(z)+\epsilon z.$
Suppose that $|f_1(z)+\epsilon z-f_2(z)| \leq |\epsilon|^3$
for all $|z|\leq \delta.$
Then, we get that $|w_1-w_2| \leq |\epsilon|^3.$ Hence, by the
previous Lemma, we have that $|f_1(z)-f_2(z)| \leq C|\epsilon|^{\frac{3}{2}}$
for all $|z| \leq \delta.$ Hence if $|z|=\delta,$

\bea
|\epsilon|^3 & \geq &  |f_2(z)+\epsilon z-f_1(z)|\\
& \geq & |\epsilon |\delta-   |f_2(z)-f_1(z)|\\
& \geq & |\epsilon| \delta- C|\epsilon|^{\frac{3}{2}}\\
& \geq & |\epsilon|(\delta-C \sqrt{|\epsilon|}),\\
& \Rightarrow & \\
\epsilon_0^2 & \geq & |\epsilon|^2 \geq  \delta-C \sqrt{|\epsilon|}
 \geq  \delta-C \sqrt{\epsilon_0}\\
& \Rightarrow & \\
\epsilon_0^2+C \sqrt{\epsilon_0} & \geq & \delta,\\
\eea

\noindent a contradiction if $\epsilon_0$ is small enough.\\

\CQFD

The following lemma is well known.

\begin{lem}
a) There is a constant $0<c<1$ so that the following holds:
Let $g$ be a holomorphic function on the unit disc
with $|g|<1$ and suppose that $g$ has $N$ zeroes in $\Delta(0,1/2).$
Then $|g|\leq c^N$ on  $\Delta(0,1/2).$\\
b) Let $g$ denote a holomorphic function on the unit disc and suppose
that $|g|<1$ and that $|g|< \eta<1$ on  $\Delta(0,1/4).$ Then
 $|g|< \sqrt{\eta}$ on  $\Delta(0,1/2).$
\end{lem}

{\bf Proof:} To prove a) set $$c= sup_{|\alpha|\leq 1/2, |z|\leq 1/2}
\frac{|z-\alpha|}{|1-z\overline{\alpha}|}<1.$$ To prove b)
observe that $\log |g|< \log \eta$ when $|z| < \frac{1}{4}.$
Hence by subharmonicity
$$
\log |g| \leq \max\{ \log \eta \frac{\log|z|}{\log 1/4}, \log \eta\}.
$$

This implies that if $|z|=1/2$, then $\log |g| \leq \frac{\log \eta}{2}.$

 \CQFD\\

{\bf Continuation of the Proof of Theorem 3.4:} Pick $\rho>0.$
Let $p\in X$. Since every leaf is dense, there is a (nonunique) continuous
curve $\gamma_p(t), 0 \leq t \leq 1$ from $\gamma(0)=p$ to a point $
\gamma_p(1)=(0,w_p)\in B_1'$ which
is contained in the leaf through $p.$ By continuity, for every
$q\in K$ close enough to $p,$ the curve $\gamma_q$  can be chosen
so that $\mbox{dist}(\gamma_q(t),\gamma_p(t))\leq \rho$
for all $0\leq t\leq 1.$

\medskip
A chain of flow boxes is a finite collection $C=\{B_{i(j)}\}_{j=1}^k$.
Let $p \in X$. We say that the leaf through $p$
follows the chain $\{B_{i(j)}\}_{j=1}^k$ if there are local
leaves $L_j \subset B_{i(j)}, \hat{L}_j:=L_j \cap B'_{i(j)},
p\in \hat{L}_1, \hat{L}_j \cap \hat{L}_{j+1}\neq \emptyset\; \forall \; j<k,
i(k)=1.$\\

By compactness there are finitely many chains of flow boxes $C_1,
\dots, C_\ell$ such that for each $p\in X$, there is an open neighborhood
$U(p)$ and a chain $C_r$ so that the leaf through $q$
follows $C_r$ for any $q\in U(p)\cap X.$\\

We will apply Lemma 3.6 repeatedly along a chain. We need to apply Lemma
3.6
at most a fixed number of times $m$ depending on the length
of each chain.
Note that every time we switch flow box there is a change of coordinates
which distorts distances by at most a factor $C>1.$

\begin{lem}
Let $\epsilon$ be sufficiently small and suppose
that $N=N(\epsilon)$ is an integer such that
$C^2c^{\frac{N}{2^m}}\leq |\epsilon|^3.$ Then
no local leaves of the laminations $L_1,\Phi_\epsilon(L_0)$
can intersect more than $N$ times in any flow box.
\end{lem}

{\bf Proof:}
Suppose that  local components of 
$L_1$ and $\Phi_\epsilon(L_0)$ intersect in more than
$N$ points in some local flow box $B'_s$. Then these local graphs
differ by at most $c^N.$ Using Lemma 3.6 they differ
by at most $c^{\frac{N}{2}}$ in a suitable larger flow box.
Changing to the coordinates of another flow box might increase
the difference to $C c^{\frac{N}{2}}.$ Applying Lemma 3.6
the difference increases to at most $C^{\frac{1}{2}}
c^{\frac{N}{4}}$ and after another change of flow box
to $C^{\frac{3}{2}}c^{\frac{N}{4}}.$ Following the leaves along a chain of flow
boxes we see inductively that the distance between continuations
of the leaves grows at most
like $C^2 c^{\frac{N}{2^k}}$ after $k$ steps.
Hence once we are in the first flow box,
the leaves differ
by at most $|\epsilon|^3.$ By the above
lemmas, this is impossible for any pair of leaves.\\

\CQFD\\

There is a constant $A$ so that for all small enough $\epsilon$
local leaves $L_1,\Phi_\epsilon(L_0)$ have at most
$N_\epsilon:=A \log \frac{1}{|\epsilon|}$ intersection points. The contruction
is stable under small perturbations by $\Psi_1,\Psi_2$
close to the identity.\\

\CQFD\\

\section{Construction from discs. Ahlfors type construction.}

In this paragraph we consider a laminated set $(X,{\mathcal L},E)$
in a compact complex manifold $M.$
 
We want to construct harmonic currents using the Ahlfors exhaustion
technique.

\subsection{When leaves are not uniformly Kobayashi hyperbolic}

We consider only the case when $X$ is not a compact
Riemann surface, possibly singular. Consider the universal covering
for each leaf. We can assume that the covering is $\CC$ or the
unit disc $\Delta.$ Let $\phi:\Delta \rightarrow L$ be a covering map.
If $|\phi'(0)|$ is not uniformly bounded, then using the Brody
technique, one can construct an image of $\CC$. The part of the image not
in $E$ is locally
contained in a leaf of $X$.

The Ahlfors exhaustion technique furnishes a positive closed current
of mass $1$ directed by the lamination. So we get the following proposition:

\begin{proposition}
Let $(X,{\mathcal L},E)$ be a laminated set with singularities. If there
is no positive closed current on $X,$ laminated on $X \setminus E$, then
there is a constant $C$ such that $|\phi'(0)| \leq C.$ 
\end{proposition}

When the $|\phi'(0)|$
are uniformly bounded, we say that the leaves are uniformly hyperbolic.

\subsection{The case with no positive closed current
directed by $(X,{\mathcal L},E)$} Let $X$ be a minimal laminated compact 
set in $\PP^2$. Suppose that
$X$ does not contain any non constant holomorphic image of $\CC$. Let $B_i$ 
be a 
covering of $K\setminus E$ by flow boxes which in local coordinates are of
the form $w_i=f_\alpha(z_i), |z_i|<1$ [but the graphs extend uniformly
to $|z_i|<2$]. Let $\Phi:\Delta\rightarrow L$ denote the universal
covering of an arbitrary leaf. We say that $x\in \Delta$ is a center point
if $\Phi(x)=(0,w_i)$ in some $B_i.$ We can normalize $\Phi$
for any center point, say $\Phi_x:\Delta\rightarrow L, \Phi_x(0)=\Phi(x).$
[i.e. we move $x$ to $0$ with an automorphism of the unit disc.]
Let $w_i=f_x(z_i)$ denote the associated graph in the flow box.
 Denote by $U_x:=\Phi_x^{-1}(\{(z_i,f_x(z_i));|z_i|<1\}).$ If $\Phi$
is a multisheeted covering, we let $U_x$ denote the connected
component containing $0.$ Then $U_x\subset
\Delta$ is a relatively compact open subset of $\Delta$ containing
$0.$ Let $0<r_x\leq R_x<1$ denote the largest, respectively smallest
radii such that $\Delta(0,r_x)\subset U_x\subset \Delta(0,R_x).$

\begin{lem}
$R_x\leq 1/2$.
\end{lem}

{\bf Proof:} Since $\Phi^{-1}_x$ maps $\Delta(0,2)$ into $\Delta(0,1)$
and sends $0$ to $0$ this follows from the Schwarz' Lemma.\\

\CQFD\\

\begin{lem}
Fix a finite number of flow boxes $\overline{B}_i \cap E=\emptyset,
i=1,\dots,\ell.$ Then $\inf \{r_x; x\; {\mbox{is a center point}},
\phi(x) \in \cup_{i=1}^\ell B_i\}>0.$
In fact the same holds if we inf over all leaves and all covers of the leaves
by discs.
\end{lem}

{\bf Proof:}
Fix an $x$ and a covering $\Phi:\Delta \rightarrow L$ of the leaf
through $x.$ Note that $\Phi_x^{-1}(\Delta(0,1)) \supset \Delta(0,r_x)$
hence by the Koebe 1/4 Theorem, $[\Phi_x^{-1}]'(0)=\alpha, |\alpha|\leq 4r_x.$
Hence  $\Phi'_x(0)=\frac{1}{\alpha}, \left| \frac{1}{\alpha}\right|
\geq \frac{1}{4r_x}.$ If $r_x \rightarrow 0$ then using the Brody
technique we construct an image of $\CC$ contained in $X.$\\

\CQFD\\

Next we prove a density Theorem for the above minimal laminations
with only Kobayashi hyperbolic leaves:

\begin{theo}
Assume $E =\emptyset.$ Fix a finite cover of $X $
by flow boxes $B_i.$
There are constants $R, N$ so that if $\Phi:\Delta\rightarrow X$
is a covering of any leaf, then $\Phi(\Delta_{kob}(x,R)\cap B_i)$
intersects at most $N$ of the graphs in $\{|z_i|<\frac{3}{2}\}$ and contains
at least one complete graph over $|z_i|<1.$
\end{theo}

{\bf Proof:} 
Let $\eta\in \cup B_j$ and $L_\eta$ the leaf through $\eta.$
There exist finitely many curves $\gamma_i^\eta\subset L_\eta,
\gamma_i^\eta(0)=\eta,
\gamma_i^\eta(1)\in B_i.$ In fact, for every $\zeta\in U(\eta)\cap X$
there are curves $\gamma_i^\zeta \subset L_\zeta$ depending continuously
on $\zeta$ and landing in $B_i$ for a small enough neighborhood $U(\eta).$
There exists a finite number $M_\eta$ so that if $\zeta\in U(\eta)\cap X,
t\in \left[\frac{j-1}{M_\eta},\frac{j}{M_\eta}\right],$ then
$\gamma_i^\zeta(t)\in \Delta_{i,\zeta,j}$ where $\Delta_{i,\zeta,j}$
is one of the unit discs in one of our chosen local flow boxes. 
It follows that if $\Phi:\Delta\rightarrow
L$ is any covering of any leaf and $\Phi(x_0)\in U(\eta),$ then
$\Phi(\Delta_{kob}(x_0,R_\eta))$ contains a disc in each flow box if $R_\eta$
is
large enough. By compactness of $X,$ there is an $R={\mbox{max}}\; R_\eta>0$
so that
if $\Phi:\Delta \rightarrow X$ is a covering of any leaf and $x \in 
\Delta$, then $\Phi(\Delta_{kob}(x,R))$ contains at least one
complete graph in each flow box. There is by Lemma 4.2 an $R'>R$
so that if $\zeta\in \Phi(\Delta_{kob}(x,R)) \cap B_i$ (or a slight extension),
then the whole graph in $B_i$ is contained in $\Phi(\Delta_{kob}(x,R')).$
It follows from Lemma 4.3 that there is an $N$ so that no
$\Phi(\Delta_{kob}(x,R))$ can intersect any flow box in more than
$N$ graphs.\\

\CQFD\\

\begin{coro}
Let $R$ be sufficiently large. Then if $D$ is any maximal disc
of center $p$ contained in an annulus $1-1/R^k<|x|<1-1/R^{k+1}$, 
and assume $\Phi(p) \in \cup B_j$, then the image
$\Phi(D)$ always will contain a full disc in each $B_i$ and will intersect
at most $N$ graphs in any $B_i$.
\end{coro}

{\bf Proof:} This follows since the Poincar\'e radius of these discs
increase to infinity independently of $k$ as $R\rightarrow \infty$.\\

\CQFD\\

\begin{coro} For large $R$ the area of $\Phi(\{
 1-1/R^k<|x|<1-1/R^{k+1}\})$ grows like $R^k$ as $k \rightarrow 
\infty.$
Moreover the length of $\Phi(\{|x|=1-1/R^{k+1}\})$ is on the
order of $R^k,$ and $|\Phi'(x)| \sim \frac{1}{1-|x|}$ near
$\partial \Delta,$ so
$$
 \int_\Delta(1-|x|)|\Phi'(x)|^2d\lambda(x)=\infty
$$
\end{coro}

{\bf Proof:} The map $\Phi$ is an isometry for the Kobayashi distance.
The Kobayashi distance  on leaves is comparable to any given Hermitian
metric. The map $\Phi$ expands a disc of radius $\epsilon$
centered at $1-2\epsilon$ to a disc of radius about $1$ in the
Kobayashi metric. Hence $|\Phi'(x)| \sim \frac{1}{1-|x|}.$\\

\CQFD\\

\begin{theo}
Let $\phi:\Delta \rightarrow M$ be a holomorphic map where $M$ 
is a compact complex Hermitian manifold. If
$$
 \int_\Delta(1-|x|)|\phi'(x)|^2d\lambda(x)=\infty,
$$
\noindent then there is a positive harmonic current $T$, supported on
$\overline{\phi(\Delta)}.$
If $\phi(\Delta)$ is contained in a leaf of a lamination ${\cal L}$, then the current
$T$ is directed by the lamination.
\end{theo}

{\bf Proof:}
Assume that $\phi(0)=p.$ Define 
$$
G_r(x):=\frac{1}{2\pi}\log^+ \frac{r}{|x|}, T_r:=(\phi)_*(G_r[\Delta]
), r<1.
$$ 
If $\theta$ is a $(1,1)$ test form on $M$
$$
<T_r,\theta>=\frac{1}{2\pi}\int_{\Delta}\log^+ \frac{r}{|x|}\phi^*(\theta).
$$
So $T_r$ is positive of bidimension $(1,1)$. The mass of $T_r$
is comparable to \\
$\int_{\Delta}\log \frac{r}{|x|}|\phi'(x)|^2 \sim
\int_{\{|x|<r\}}(r-|x|)|\phi'(x)|^2=:A(r).$ A direct computation gives
$$
i \partial \overline{\partial} T_r=\phi_*(\nu_r)-\delta_p,
$$
\noindent where $\delta_p$ is the Dirac mass at $p$ and $\nu_r$ is the Lebesgue measure on the circle
of radius $r.$ Let $T$ be a cluster point of $T'_r:=\frac{T_r}{A(r)}.$
Since $A(r) \rightarrow \infty$ we have $i \partial \overline{\partial}T=0.$

\medskip

If $\phi(\Delta)$ is a leaf of a lamination ${\cal L}$, there is a $(1,0)$
form $\gamma$ such that $[\phi(\Delta)]\wedge \gamma=0.$ Hence
$T_r \wedge \gamma=0.$ 
Hence $<T,i\partial_b\overline{\partial}_bf>=0$ and we have a decomposition
of $T$ as in Theorem 2.18.\\

\CQFD\\

\begin{remark} If we assume that $\lim_{r \rightarrow 1}
(1-r) \int_{D_r} |\phi'(z)|^2=\infty,$ it follows from a result of Ahlfors
that $\underline{\lim}_{r\rightarrow 1} \frac{\ell(\phi(\partial
\Delta_r))}{{\mbox{Area}}(\phi(\Delta_r))}\rightarrow 0,$ $\ell$ represents
length. Hence one can choose the current $T$ to be closed.\\
In the usual Ahlfors procedure to construct a closed current
starting from an image of $\CC$ one has to first extract good subsequences from
$\frac{\Phi_*[\Delta_R]}{{\mbox{Area}}\; \Phi_* [\Delta_R]}$ when
$R \rightarrow \infty.$ Then cluster points of these give closed
currents. In our case, there is no need to first take a subsequence.
\end{remark}

\begin{theo}
There is no nonvanishing holomorphic 
vector field
along leaves on a laminated
compact with only hyperbolic leaves.
\end{theo}

{\bf Proof:} If we pull back the vector field to a disc covering
a leaf, we get a holomorphic function on the unit disc going uniformly
to infinity at the boundary, since as we have seen $|\phi'(x)| \sim \frac{1}
{1-|x|}$ and the vector field is bounded in flow boxes.\\

\CQFD\\

\begin{proposition}
If a positive harmonic current gives mass to a leaf, then this leaf is
a compact Riemann surface.
\end{proposition}

\begin{lem}
Let $T$ be a laminated harmonic current. Let $\phi$
denote the covering map $\Delta \rightarrow L.$ If $H$ denotes the
analytic continuation of $h \circ \phi$ in a flow box, then we have the
estimate:
$$
c (1-|x|) \leq H  \leq C \frac{1}{(1-|x|)}.
$$
\end{lem}

{\bf Proof:}
The estimate 
follows from Harnack's inequality and the Hopf Lemma.\\

\CQFD\\

{\bf Proof of the Proposition:} We assume first that the leaf is hyperbolic.
Let $\phi, H$ be as in the Lemma.
Suppose at first that $H$ is unbounded. Then we can choose a sequence
$p_n \rightarrow \partial \Delta$  and $H(p_n) \rightarrow \infty$ such
that $H$ is uniformly large on $\Delta(p_n,R)$
by Harnack, where $R$ is as in Theorem 4.4.
Hence $T$ will have infinite mass on a flow box.
If $H$ is bounded and nonconstant, we can choose $\theta_n$ so that
$lim_{r \rightarrow 1}H(re^{i \theta_n})$ are different.
We can again choose $p_n \rightarrow \partial \Delta$  so that
$\phi(\Delta(p_n,R))$ are disjoint and again $T$ will have infinite mass
on a flow box. If $H$ is constant, we get a positive
closed current and the leaf has an analytic closure. The same
argument applies to the case when the leaf is not hyperbolic.\\

\CQFD\\

\section{Vanishing of $\int T \wedge T.$}

\begin{definition}
The harmonic laminated current has finite transverse energy if
in some local flow box
$\int \log |\alpha-\beta| d\mu(\alpha) d\mu(\beta)<\infty.$
\end{definition}

Having finite transverse energy is well defined and independent of the 
choice of flow box.

\medskip

Recall that $\Phi_\epsilon([z:w:t])=[z:w+\epsilon z:t].$ If $T$
is a current, let $T_\epsilon:= (\Phi_\epsilon)_*(T).$

\begin{theo}
If a harmonic current for a
laminated compact in $\PP^2$ has finite transverse energy, then
the geometric intersection $T \wedge T_\epsilon \rightarrow 0.$
The same conclusion holds for ${\mathcal C}^1-$laminations without the
hypothesis of finite transverse energy. In both cases we have
$\int T \wedge T=0.$  
\end{theo}

{\bf Proof:}
We calculate the geometric wedge product $T \wedge T^\epsilon$ in a flow box.
Set $T=\int h_\alpha[\Delta_\alpha]d\mu(\alpha),
T^\epsilon =\int h_\beta^\epsilon[\Delta^\epsilon_\beta]d\mu(\beta).$
Let $\phi$ be a test function supported in a flow box.
To avoid confusion, we index with $g$ when wedge products are geometric
during the proof. 
We have

$$
<T \wedge T^\epsilon,\phi>_g=
\int \sum_{p \in J_{\alpha,\beta}}\phi h_\alpha(p)h^\epsilon_\beta(p)
d\mu(\alpha)d\mu(\beta)$$

\noindent where $J_{\alpha,\beta}$
consists of the intersection points of $\Delta_\alpha$ and
$\Delta^\epsilon_\beta.$ Assume at first that $\mu$ has finite
transverse energy. Using the estimate on the size of $J_{\alpha,\beta}$
in Theorem 3.3, we get:

\bea
|(T \wedge T^\epsilon)_g (\phi)| & \leq & C_1
\|\phi\|_\infty   \int_{{\mbox{dist}}(\Delta_\alpha,\Delta_\beta)\leq
C \epsilon} A \log \frac{1}{|\epsilon|}d\mu(\alpha)d\mu (\beta)\\
& \leq & C_2\|\phi\|_\infty
\int_{{\mbox{dist}}(\Delta_\alpha,\Delta_\beta)\leq C \epsilon}
 \log \frac{1}{dist(\Delta_\alpha,\Delta_\beta)}d\mu(\alpha)d\mu(\beta)\\
& \rightarrow & 0 \; as\; \epsilon \rightarrow 0\\
\eea

In the ${\mathcal C}^1$ case the number of intersection points is bounded by 
$N$ independent of $\epsilon$ (Theorem 3.1). Hence

\bea
|(T \wedge T^\epsilon)_g(\phi)| & \leq &
C \|\phi\|_\infty \int_{{\mbox{dist}}(\Delta_\alpha,\Delta_\beta)
\leq C\epsilon}Nd\mu(\alpha)d\mu(\beta)\\
& \rightarrow & 0\\
\eea

\noindent since $\mu$ has no pointmasses by Proposition 4.10.
Next we show that $Q(T,T)=\int T \wedge T=0.$ It suffices to show
by Theorem 2.12 that $Q(T,T_\epsilon) \rightarrow 0$ or even that
for smoothings $T^\delta, T_\epsilon^{\delta'}$ that
$Q(T^\delta,T_\epsilon^{\delta'}) \rightarrow 0$ when
$\delta,\delta'$ are sufficiently small compared to $\epsilon$
and $\delta,\delta',\epsilon \rightarrow 0.$ 

\medskip

Note that the estimate on the geometric wedge product is stable under
small translations of $T,T_\epsilon$. This is what allows us to
smooth.

Let $\phi$ be a test function supported in some local
flow box. 
As above, the value of the geometric wedge product
on $\phi$ is:

$$
<T \wedge T^\epsilon,\phi>_g=
\int \sum_{p \in J_{\alpha,\beta}}\phi h_\alpha(p)h^\epsilon_\beta(p)
d\mu(\alpha)d\mu(\beta)$$

We can write this as

$$
<T \wedge T^\epsilon,\phi>_g=
\int \left(
\int_{\Delta^\epsilon_\beta}[\phi h_\alpha h^\epsilon_\beta](p)
i\partial \overline{\partial} \log|w-f_\alpha(z)|\right) 
d\mu(\alpha)d\mu(\beta)$$

The same applies when we do this for slight translations within small
neighborhoods $U(\epsilon)$ of the identity in $U(3)$
and
their smooth averages $T^\delta,$

$$
<T^\delta \wedge T_\epsilon,\phi>_g=
\int \left(
\int_{\Delta^\epsilon_\beta}[\phi h^\epsilon_\beta](p)
T^\delta \right) d\mu(\beta)=<T_\epsilon,\phi T^\delta>.$$

Averaging also over small translations of $T^\epsilon$ we get

$$
<T^\delta \wedge T^{\delta'}_\epsilon,\phi>_g=
<T^{\delta'}_\epsilon,\phi T^\delta>.$$

We still have that $<T_\epsilon^{\delta'},\phi T^\delta> \rightarrow 0$
when $\delta,\delta' << \epsilon, \epsilon \rightarrow 0.$ If we apply this
to $\phi=1$, we get
$<T^{\delta'}_\epsilon,T^\delta>=Q(T^{\delta'}_\epsilon,T^\delta)
\rightarrow 0.$ Hence
$Q(T,T)=0.$\\

\CQFD\\

\begin{coro}
If a laminated compact set in $\PP^2$
carries a positive closed laminar 
current $T$, then $T$ has infinite transverse energy.
\end{coro}

{\bf Proof:}
If $T \neq 0$ has finite energy, then $0=\int T \wedge T=
|\int T \wedge \omega|^2-E(T,T)$ but $E(T,T)=0$
since $T$ is closed. Hence $T=0$, a contradiction.\\

\CQFD\\

J. Duval has independently obtained this Corollary.
Hurder and Mitsumatsu proved that there is no ${\mathcal C}^1$
lamination in $\PP^2$ which carries a positive closed current
[HM].

\begin{coro}
If $(X,{\mathcal L})$ is a ${\mathcal C}^1$ lamination on $\PP^2$
with only hyperbolic leaves,
then $$T=\lim_{r \nearrow 1}\phi_*\frac{
\left(\log^+ \left(\frac{r}{|z|}\right)[D_r]\right)}{A_r}$$
\noindent uniformly with respect to $\phi.$
\end{coro}

{\bf Proof:}
We know from [HM] that there is no positive closed current directed by 
${\mathcal L}$. It follows from Corollary 2.23 and Theorem 5.2 that there
is a unique harmonique current of mass $1$ on $(X,{\mathcal L})$.
Hence the result follows.\\

\CQFD\\

\section{Examples of harmonic current}

In this section we investigate harmonic currents on $\PP^2$
of the form $T=i \partial u \wedge \overline{\partial u}.$
Our main result is that if $u \in {\cal C}^2_{\RR}(\PP^2)$
and $i \partial \overline{\partial} T=0,$ then $u$ is constant,
hence $T\equiv 0.$ We also compute the energy of some positive
harmonic currents. 

\medskip

Let $M$ be a complex manifold of dimension $m.$
For $1 \leq k \leq m,$
 we define ${\cal P}^{(k)}_-(M)$ as the cone of upper semicontinuous
real functions $v$ on $M$ such that for every $p\in M,$
there is an open neighborhood $U$ of $p$ and $\{v_n\}\subset {\cal C}^2(U)$ such
that $v_n \searrow v$ in $U$ and
$(-1)^k(i\partial \overline{\partial} v_n)^k \leq \epsilon_n \omega^k, \epsilon_n \searrow 0.$ We say
that $U$ is associated to  $v.$ Here 
$\omega$ denotes a strictly positive hermitian form. 
Notice that this condition implies that when $\epsilon_n=0,$
not all eigenvalues of $i\partial \overline{\partial} v_n$ can be strictly negative.
We define ${\cal P}^{(k)}_+:=-{\cal P}^{(k)}_-(M)$ and ${\cal P}^{(k)}
={\cal P}_+^{(k)}\cap {\cal P}_-^{(k)}.$
In particular, $P^{(1)}_-$ consists of the plurisubharmonic functions
and $P^{(1)}$ are the pluriharmonic functions. In dimension $2$, a
smooth function
$v$ belongs to $P^{(2)}_-$ if and only if its Levi form has at most one
eigenvalue of each sign. These functions then, also belong to 
$P^{(2)}_+$ and hence $P^{(2)}.$ Pseudoconvex domains are usually 
characterized by plurisubharmonic functions, i.e. $P^{(1)}_-.$ 
We show here that $P^{(2)}_-$ works as well, and that there are
similar results for $P^{(k)}_-, k>2.$

 Let $\Phi:M \rightarrow N$ be a holomorphic map between complex manifolds.
If $v\in {\cal P}^{(k)}_-(N),$ then
$v \circ \Phi \in {\cal P}_-^{(k)}(M).$ In particular, if $k> {\mbox{dim}}\;
N,$ any upper semicontinuous $v$ is in $P_-^{(k)}(N),$ hence
$v \circ \Phi \in P^{(k)}_-(M).$

\medskip

We give some examples of compact complex 
manifolds for which $P^{(2)}(M) \neq \RR:$\\

\smallskip

1. Tori: Let $T$ be a torus. Then $T=\CC^k$ mod a lattice
generated by $\{v_i\}_{i=1}^{2k}.$ Let $\pi:\CC^k \rightarrow \RR$,
$\pi(\sum x_i v_i)=x_i.$ Let $v=\phi(x_1)$ where $\phi$ is a smooth function supported in
$]0,1[.$ Then $(i\partial \overline{\partial} v)^2=0.$\\

\smallskip

2. Hopf manifolds. Let $M=\CC^2/<\phi>$ where $<\phi>$ denotes the
group generated by $\phi(z_1,z_2)=(\alpha_1 z_1,\alpha_2 z_2)$ with
$\alpha_1,\alpha_2$ fixed, $0<|\alpha_1| \leq |\alpha_2| <1.$ Fix
$r$ such that $|\alpha_1|=|\alpha_2|^r.$ Define
$$v(z_1,z_2)=\frac{|z_1|^2}{|z_1|^2+|z_2|^r}$$
The function $v$ is well defined on $M$ and $(i\partial \overline{\partial} v)^2=0.$\\

\smallskip

3. For any surface admitting a projection on a Riemann surface $P^{(2)}(M)
\neq \RR,$ for example ruled surfaces. Actually the Hopf surfaces above admit
such a projection $\CC^2/(\phi) \rightarrow \PP^1,$ 
$(z_1,z_2) \rightarrow [z_1^q:z_2^p]$ if
$\alpha_1^q=\alpha_2^p.$\\

\medskip

For $k>1,$ set $z=(z_1,\dots,z_{k-1}), |z|= \max \{|z_1|,\dots,|z_{k-1}|\}$
and let $H^r_{k-1}$ denote the Hartogs figure:

$$
H^r_{k-1}:=\{(z,w)\in \CC^{k-1} \times \CC= 
\CC^k; |z|\leq 1+r,|w|\leq 1\}\setminus \{r<|w|\leq 1, |z|< 1\}.
$$

Let $\hat{H}^r_{k-1}:= \{(z,w)\in \CC^k; |z|\leq 1+r,|w|\leq 1\}$.

\begin{definition} Let $2 \leq k \leq m.$ We say that an open set 
$N\subset M$ is\\
$(k-1)-$pseudoconvex if whenever  $\Phi:U \rightarrow M$ is a biholomorphic 
map
of a neighborhood $U\supset \hat{H}^r_{k-1}$ onto its image and
$\Phi(H^r_{k-1}) \subset N,$ then $\Phi(\hat{H}^r_{k-1}) \subset N.$
\end{definition}

\begin{remark} In the case $k=2$, the definition is  
equivalent to $N$ being pseudoconvex.
\end{remark}

\begin{prop}
Let $M$ be a complex manifold of dimension $m \geq 2.$ Let $N$ be a
connected open 
set in $M$. Assume $v\in P_-^{k}(M), 2\leq k \leq m$ and $v<0$
on $N,$ $v_{|\partial N}\equiv 0, v \leq 0$ on $M.$
Let ${\cal U}=\{U_\alpha\}$ be a cover of $M$ associated to $v.$ If
$\Phi(\hat{H}^r_{k-1})\subset U_\alpha$ and $\Phi({H}^r_{k-1})\subset N,$
then $\Phi(\hat{H}^r_{k-1})\subset N.$ In particular if $k=2$ then $N$
is pseudoconvex in $M$.
\end{prop}

{\bf Proof:} Assume that $\Phi({\hat{H}}^r_{k-1})\setminus N
\neq \emptyset.$  Then also, $\Phi({\hat{H}}^r_{k-1})\cap
\partial N \neq \emptyset.$ There is a biholomorphic map
$\Phi:U \rightarrow M$ onto its image, ${\hat{H}}^r_{k-1}\subset U,$
$\Phi(H^r_{k-1}) \subset N.$
We can then assume that $M=U,$ $\Phi={\mbox{Id}},$
$H^r_{k-1}\subset N$ and that there is an interior point of ${\hat{H}}^r_{k-1}$
in $\partial N.$ 
Assume that $v_n \searrow v,v_n \in {\cal C}^2(U), (-1)^k(i\partial \overline{\partial} v_n)^k \leq 
\epsilon_n \omega^k.$
We also assume that $\partial N$ contains a point $(z_0,w_0),
|z_0|,|w_0|<1.$ 
Then $v<0$ on ${H^r_{k-1}} $ and $v(z_0,w_0)=0.$ This is where we use
that $v_{|\partial N} \equiv 0.$ We would like
to get a contradiction.

Let $X:=\{|z|\leq 1, r< |w|\leq 1\}.$
Define the function $u$ by

$$
u(z,w):= \eta(1-\sum_{j=1}^{k-1}|z_j|^2)-\frac{\epsilon}{|w|^2}+\delta
$$

\noi where $\eta<<\epsilon<\delta<<1.$ We will choose the constants so that
$v < u$ on $\partial X$ and $u(z_0,w_0)<0.$ First observe that if $\delta$ is small enough,
then automatically $v<u$ on all of $\partial X$ except possibly where
$|w|=1$ and $|z|<1.$ Fix any such $\delta.$ Let $\epsilon<\delta$ be chosen
big enough that $-\frac{\epsilon}{|w_0|^2}+\delta<0.$ Since
$\epsilon<\delta$ and $v\leq 0$ on $|w|=1$ and $|z|<1$,
 $v<u$ on all of $\partial X$ if we choose $\eta=0$.
To finish the choice of constants, $\eta, \epsilon,\delta$ choose
$\eta>0$ small enough that $v<u$ still on $\partial X$ and in, addition 
$u(z_0,w_0)<0.$ Then $i\partial \overline{\partial} (-u) \geq a \omega$ on $X$ for some constant
$a>0.$ 

\medskip

Next choose $n$ large
enough so that $v_n<u$ on $\partial X$. We have $v_n(z_0,w_0)
\geq v(z_0,w_0)=0>u(z_0,w_0).$ Then if we add a strictly
positive constant $c_n$ to $u$ we can assume that

\bea
v_n & < & u+c_n \; {\mbox{on}}\; \partial X\\
v_n & \leq & u+c_n \; {\mbox{on}}\; X\\
v_n(z_1^n,w_1^n) & = & u(z_1^n,w_1^n)+c_n, (z_1^n,w_1^n) \in X\\
\eea

\bigskip

This implies that $i\partial \overline{\partial} v_n(z_1^n,w_1^n)\leq i 
\partial \overline{\partial} u(z_1^n,w_1^n).$ 
Hence $i\partial\overline{\partial}  (-v_n)(z_1^n,w_1^n) \geq a\omega.$ This implies
that

\bea
\epsilon_n \omega^k & \geq & (-1)^k (i\partial \overline{\partial} v_n)^k(z_1^n,w_1^n)\\
& = &  (i \partial \overline{\partial} (-v_n))^k(z_1^n,w_1^n)\\
& \geq & a^k \omega^k,\\
\eea

\noindent a contradiction.\\

\CQFD\\

\begin{coro}
If $M$ is a compact manifold and $v\in P^{(2)}(M),$ then
$K(v): =\{p;v(p)=\max_M v\}$
is pseudoconcave, i.e. $M \setminus K(v)$ is
pseudoconvex. If $v \in P^{(k)}(M) \cap {\cal C}^2,$ then
$M\setminus K(v)$ is $k-1$ pseudoconvex.
\end{coro}

\begin{remark}
Let $v$ be a ${\mathcal C}^2$ function on a compact complex manifold
of dimension $m.$ Stokes' Theorem implies that if $(-1)^k
(i \partial \overline{\partial} v)^k \leq 0,$ then $(dd^c v)^k =0.$ In the case
of compact K\"{a}hler manifolds, Stokes' Theorem applied to
$(-1)^k(i\partial \overline{\partial} v)^k \wedge \omega^{m-k}$ shows that the same
conclusion holds.
\end{remark}

\begin{remark}
The proof above shows that if an upper semicontinuous function $v$
is locally a decreasing limit of ${\cal C}^2$ functions $v_n$ such that
at each point $i\partial \overline{\partial} v_n$ has $m-1$ nonnegative eigenvalues,
then $K(v)$ is pseudoconcave. Namely, we get by the above construction
with a Hartogs figure of dimension two:
\bea
v_n & < & u+c\; {\mbox{on}}\; \partial X,\\
v_n & \leq & u+c  \; {\mbox{on}}\;X\\
v_n(z_1,w_1) & = & u(z_1,w_1)+c \; {\mbox{at some point of }}\; X.\; 
{\mbox{Hence}}\\
i \partial \overline{\partial}(v_n-u)(z_1,w_1) & \leq & 0\\
\eea
\noi which contradicts that one eigenvalue is nonnegative.
\end{remark}

\begin{coro}
Let $v$ be a continuous function on $\PP^m$ such that $v\in {\cal P}_-^{(2)}$,
then $v$ is constant. In particular there are no
nonconstant functions in ${\mathcal C}^2_{\RR}$ such that $T=
i\partial v \wedge \overline{\partial v}$ satisfies
$i \partial \overline{\partial} T=0.$
\end{coro}

{\bf Proof:}
We know that $K(v)$ is pseudoconcave. We show next that also
$K(-v)$ is pseudoconcave. It suffices to show that $-v$ is also
a decreasing limit of ${\cal C}^2$ functions $w_n$, 
$(i\partial \overline{\partial} w_n)^2 \leq \epsilon_n \omega^2.$
For this, let $v_n$ be such a sequence for $v.$ Taking a subsequence if
necessary we can assume that $v \leq v_n \leq v+1/2^n.$ Set $w_n=-v_n+1/n.$
 Since on $\PP^m$ the Levi problem has a positive
solution, this implies that the complements of $K(v)$ and $K(-v)$ are
both Stein. But then the intersection of the two domains is a Stein
manifold of dimension $>1$ with two ends, unless $K(v)$ and $K(-v)$ have
a nonempty intersection. But then $v$ must be constant.\\

\CQFD\\

\begin{coro}
If the Levi problem is solvable on a compact complex manifold $M$, then
$P^{(2)}$ only contains constant functions.
Hence there is no nonconstant holomorphic
map from $M$ to a manifold with nontrivial
$P^{(2)}(M).$
\end{coro}

 Recall [BS] that given a positive closed current 
$S$ on $M,$ an upper semicontinuous 
function $\phi$ defined on ${\mbox{Supp}}(S)$ is $S$-plurisubharmonic if for
every $p \in {\mbox{Supp}}(S)$ there is an open set $U,$ $p\in U$ and a 
sequence  $\phi_n \in {\cal C}^2(U),$ such that
$\phi=\lim_{\searrow}\phi_n$ on $U \cap S$ and
$i\partial \overline{\partial} \phi_n\wedge S \geq 0$.
A function $\phi$ is $S-$ pluriharmonic if both $\phi$ and
$-\phi$ are \\
$S-$plurisubharmonic.

\begin{theo}
Let $S$ be a positive closed current of bidegree $(1,1)$ in $\PP^2.$ Assume
$\phi$ is $S-$plurisubharmonic and let $K(\phi)=
\{p\in {\mbox{Supp}}(S); \phi(p)= \max \phi\}. $ Then $\PP^2 \setminus
K(\phi)$ is pseudoconvex. If $\phi$ is $S-$pluriharmonic, then $\phi$ is constant.
\end{theo}

{\bf Proof:}
Recall that for $S-$plurisubharmonic functions, the local maximum principle
is valid [BS], Prop. 3.1. We claim that $\PP^2
\setminus K(\phi)$ is pseudoconvex. 
We modify the proof of Proposition 6.3. Assume that 
$\PP^2
\setminus K(\phi)$ is not pseudoconvex.
We can assume in local coordinates that $K(\phi)$
contains a point $(z_0,w_0), |z_0|,|w_0|<1$ and that $K(\phi)$
does not intersect the Hartogs figure $H=\{(z,w); |w| \leq r<1, |z|\leq 
1+\delta,\}\cup \{(z,w); 1\leq |z| \leq 1+\delta, |w|\leq 1\}.$
We can also assume that on a fixed neighborhood of $\hat H=
\{(z,w);|z|\leq 1+\delta, |w|\leq 1\}$
there is a sequence of $\mathcal C^2$ functions $\phi_n
\searrow \phi$ on Supp$(S)$, $i\partial \overline{\partial}
 \phi_n \wedge S\geq 0.$
We can assume $\phi=0$ on $K(\phi).$
Then $\phi<0$ on ${H} \cap {\mbox{Supp}}(S)$ 
and $\phi(z_0,w_0)=0.$ 

Let $X:=\{|z|\leq 1, r< |w|\leq 1\}.$
Define the function $u$ by

$$
u(z,w):= \eta(1-|z|^2)-\frac{\epsilon}{|w|^2}+\delta
$$

\noi where $\eta<<\epsilon<\delta<<1.$ We will choose the constants so that
$\phi < u$ on $\partial X \cap {\mbox{Supp}}(S)$.
First observe that if $\delta$ is small enough,
then automatically $\phi<u$ on all of $\partial X \cap {\mbox{Supp}}(S)$ 
except possibly where
$|w|=1$ and $|z|<1.$ Fix any such $\delta.$ Let $\epsilon<\delta$ be chosen
big enough that $-\frac{\epsilon}{|w_0|^2}+\delta<0.$ Since
$\epsilon<\delta$ and $\phi\leq 0$ on $|w|=1$ and $|z|<1$,
$(z,w)\in {\mbox{Supp}}(S),$ $\phi<u$ on all of 
$\partial X \cap {\mbox{Supp}}(S)$ if we choose $\eta=0$.
To finish the choice of constants, $\eta, \epsilon,\delta$ choose
$\eta>0$ small enough that $\phi<u$ still on $\partial X \cap
{\mbox{Supp}}(S)$ and in, addition 
$u(z_0,w_0)<0.$ 

\medskip

Next choose $n$ large
enough so that $\phi_n<u$ on $\partial X \cap {\mbox{Supp}}(S)$. 
Then if we add a strictly
positive constant $c$ to $u$ we can assume that

\bea
\phi_n & < & u+c \; {\mbox{on}}\; \partial X \cap {\mbox{Supp}}(S)\\
\phi_n & \leq & u+c \; {\mbox{on}}\; X \cap {\mbox{Supp}}(S)\\
\phi_n(z_1,w_1) & = & u(z_1,w_1)+c, (z_1,w_1) \in X \cap  {\mbox{Supp}}(S)\\
\eea

\bigskip

Now, $-u$ is plurisubharmonic, so $i\partial \overline{\partial}
(-u) \wedge S\geq 0.$ Hence
$\phi_n-u$ is $S-$plurisubharmonic so this contradicts
the local maximum modulus principle for $S-$ plurisubharmonic
functions.

\medskip

If $\phi$ is $S-$ pluriharmonic, then $K(\phi)$ and $K(-\phi)$ 
intersect, hence $\phi$ is constant.\\

\CQFD\\

\begin{prop}
If $v \in {\cal P}_-^{(k)}(M)$ then $v$ satisfies the local maximum 
principle.
\end{prop}

{\bf Proof:}
Recall that the local maximum principle says that
for every ball $\max_B v \leq \max_{\partial B} v.$
This follows, since the Hartogs figure argument is local. In fact, let
$K$ denote the compact set at which the maximum is reached.
Let $p\in \partial K$ and use a Hartogs figure there.\\

\CQFD\\

There are positive closed currents $T$ on $\PP^2$
of the form $T=i\partial u \wedge \overline{\partial u},$
$u$ continuous except at one point and such that $\int T \wedge
T \neq 0,$ for example: $u=\log^+|z|$ in $\CC^2$ if
$[z:w:t]$ are the homogeneous coordinates in $\PP^2.$ 

\bigskip

\begin{prop} Consider $C=\{T \geq 0,
i\partial \overline{\partial} T=0, \int T \wedge \omega=1\}.$
Then $\inf_ {T\in C} \int T \wedge T \leq 1-\frac{1}{2\pi^2}.$
\end{prop}

{\bf Proof:}
Let $u(|z|^2), v(|w|^2)$ be ${\mathcal C}^\infty$ real valued functions
with support in the unit interval. Define $\psi(z,w):=u(|z|^2)+
iv(|w|^2).$ Let $T:=i \partial \psi \wedge \overline{\partial \psi}$
on ${\mathbb C}^2.$ Then $T \geq 0$ and $T \wedge T=0.$ Moreover, $T$
is pluriharmonic on ${\mathbb C}^2.$

We want to decompose $T$ as in Proposition 2.6.

\bea
i\partial \psi \wedge \overline{\partial \psi}
& = & i(u' \overline{z}dz+iv' \overline{w}dw) \wedge
(u'zd \overline{z}-i v' wd\overline{w})\\
& = & i (u')^2 z\overline{z}dz \wedge d\overline{z}
+i (v')^2 w\overline{w}dw \wedge d\overline{w}\\
& + &  u'v' z\overline{w} d\overline{z}\wedge dw
+ u'v' \overline{z}w dz \wedge  d\overline{w}\\
\eea

Let 
\bea
U(z) & := & \frac{i}{\pi} \int \log|z-x|(u'(x \overline{x}))^2x \overline{x}
dx \wedge d \overline{x}\\
V(w) & := & \frac{i}{\pi} \int \log|w-y|(v'(y \overline{y}))^2y \overline{y}
dy \wedge d \overline{y}\\ 
\eea

Then

\bea
i \partial \psi \wedge \overline{\partial \psi}
& = & i \partial \overline{\partial} U(z)+ i \partial
 \overline{\partial} V(w)+
\overline{\partial}(uv'\overline{w}dw)-uv'' w\overline{w}
d\overline{w}\wedge dw\\
& - & uv' d\overline{w}\wedge dw
+{\partial}(uv'wd\overline{w})-uv" w\overline{w}
dw \wedge d\overline{w}-uv' dw \wedge d\overline{w}\\
\eea

Hence:

\begin{lem}
On ${\mathbb C}^2,$ 
\bea
T & = & i \partial \psi \wedge \overline{\partial \psi}\\
& = & i \partial \overline{\partial}U(z)+ i \partial \overline{\partial}V(w)
+\overline{\partial}(uv' \overline{w}dw)+\partial(uv' wd\overline{w}).\\
\eea
\end{lem}

Let $A:= i/\pi \int (u')^2|z|^2 dz \wedge d \overline{z},
B:= i/\pi \int (v')^2|w|^2 dw \wedge d \overline{w}.$
Then $U(z)=A \log |z|, |z|>1, V(w)=B \log |w|, |w|>1.$
We decompose $T$ further:
Let $h:= U(z)+V(w)-\frac{1}{2}(A+B) \log (1+|z|^2+|w|^2).$
Then

\begin{lem}
On ${\mathbb C}^2:$ 
\bea
T & = & i \partial \overline{\partial}h(z,w)+\frac{1}{2}(A+B) \omega
+ \overline{\partial}(uv'\overline{w}dw)
+\partial(uv'wd\overline{w})\\
\omega & := & i \partial \overline{\partial} \log(1+|z|^2+|w|^2).\\
\eea
\end{lem}
We extend $T$ to ${\mathbb P}^2$ as $\tilde{T}$, the trivial
extension. We need to know that $T$ has finite mass near
$[0:1:0]$ and $[1:0:0]$ for $\tilde{T}$ to be well defined.
To extend first across the line at infinity, $\eta=0$ away from
$[0:1:0]$ and $[1:0:0]$, we extend the three parts individually.
First $\omega$ extends as the K\"{a}hler form, also called $\omega.$
The form $uv'wd\overline{w}$ has compact support and extends
as $S$ trivially.
Next we investigate $h$ near $[0:1:0].$ We calculate in local coordinates.
$[z:w:1]=[Z:1:t], $ to get
$h(z,w)=\tilde{h}(Z,t))= U(Z/t)+V(1/t)-\frac{1}{2}(A+B)
\log (1+|Z/t|^2+|1/t|^2).$ When $|Z/t|>1$ we have
$U(Z/t)=A \log |Z/t|, V(1/t)=B \log |1/t|.$ So $\tilde{h}(Z,t)
=A \log|Z|-A \log |t|-B\log |t|
-\frac{1}{2}(A+B) \log (1+|Z|^2+|t|^2)+ \frac{1}{2}(A+B)\log |t|^2,$
$\tilde{h}(Z,t)=A \log |Z|-\frac{1}{2}\log(1+|t|^2+|Z|^2).$
Hence $\tilde{h}(Z,t)$ extends smoothly across $\eta=0$ except possibly
at $[0:1:0]$ and $[1:0:0].$ In particular $i \partial \overline{\partial}
\tilde{h} $ extends trivially at $\eta=0$ except
at $[1:0:0]$ and $[0:1:0].$ Next we calculate in a neighborhood
of $[0:1:0].$ We get

$$\tilde{h}(Z,t)=U(Z/t)+A \log |t|-\frac{A+B}{2}\log (1+|t|^2+|Z|^2).$$

The function $U(Z/t)+A\log|t|=:\phi(Z,t)$ is plurisubharmonic
when $t \neq 0$ and equals $A\log |Z|$ when $|Z/t|>1$ or $t=0, Z \neq 0.$
So $\phi(Z,t)$ is plurisubharmonic away from the origin. Hence $\phi$ has
a well defined plurisubharmonic extension through $(0,0)$ by setting
$\phi(0,0)=-\infty.$
It follows that $\tilde{h}$ is a global quasiplurisubharmonic function on
$\mathbb P^2$ with poles at $[0:1:0], [1:0:0].$ Hence, 

\begin{lem}
The trivial extension $\tilde{T}$ is given by
$\tilde{T}=\frac{A+B}{2}\omega+\partial S+\overline{\partial S}
+i \partial \overline{\partial}\tilde{h}, S=uv'wd\overline{w}.$
\end{lem}

It is easy to check that $\tilde{T}$ is pluriharmonic on $\PP^2.$

\medskip

{\bf End of Proof of Proposition 6.11:} This follows since $\partial \overline{\partial}\tilde{h}$
has no mass on the line at infinity. Hence $\int \overline{\partial}S
\wedge \partial \overline{S}=AB,$$
\int T \wedge T= \left(\frac{A+B}{2}\right)^24\pi^2-2AB$ and
$\int T \wedge \omega=\frac{A+B}{2} 2\pi$ so if we let
$T_1=\frac{T}{\int T \wedge \omega}$ we find
$$\int T_1 \wedge T_1 = 1-\frac{2AB}{\left(\frac{A+B}{2}\right)
4\pi^2}.$$ The minimum is reached for $A=B$ and equals
$1-\frac{1}{2\pi^2}.$\\ 

\CQFD\\

\begin{proposition}
Let $M$ be a complex surface and $\rho\in {\mathcal C}^2_{\RR}(M).$
Assume that $\partial \rho$ is non vanishing on $X=
\{\rho=0\}$ and that
$(i\partial \overline{\partial} \rho)^2=0$ on $X$ and
also that $i\partial \overline{\partial} \rho \wedge \partial \rho
\wedge \overline{\partial} \rho={\mathcal O}(\rho^2)$. Then
$T=i\delta_{\{\rho=0\}} \partial \rho \wedge \overline{\partial}\rho 
$ is a smooth positive harmonic current. Moreover $T \wedge T=0.$ 
\end{proposition}

{\bf Proof:}
Choose a $\chi\in {\mathcal C}^\infty_0(-1,1)$, $\chi \geq 0, \int \chi=1.$
Let $$T_\epsilon =i \frac{1}{\epsilon} \chi\left(\frac{\rho}{\epsilon}\right)
\partial \rho \wedge \overline{\partial} \rho.$$
Then we have 
$$
i \partial \overline{\partial} T_\epsilon
=-\frac{1}{\epsilon^2}\chi'\left(\frac{\rho}{\epsilon}\right)
\partial \rho \wedge \overline{\partial}\partial \rho \wedge
\overline{\partial} \rho-\chi \left(\frac{\rho}{\epsilon}\right)
\overline{\partial} \partial \rho 
\wedge \partial \overline{\partial} \rho.
$$

Clearly then $i \partial \overline{\partial} T_\epsilon \rightarrow 0.$\\

\CQFD\\

\bigskip

\noindent John Erik Forn\ae ss\\
Mathematics Department\\
The University of Michigan\\
East Hall, Ann Arbor, MI 48109\\
USA\\
fornaess@umich.edu\\

\noindent Nessim Sibony\\
CNRS UMR8628\\
Mathematics Department\\
Universit\'e Paris-Sud\\
Batiment 425\\
Orsay Cedex\\
France\\
nessim.sibony@math.u-psud.fr\\
\end{document}